\documentclass[twoside,12pt]{article}
\usepackage{amssymb}

\usepackage{makeidx}
\usepackage{a4wide}
\usepackage{amsthm}
\usepackage{amsmath}
\usepackage[all]{xy}
\usepackage{enumerate}
\usepackage{fancyhdr}

\newtheorem{theorem}{Theorem}[section]
\newtheorem{proposition}[theorem]{Proposition}
\newtheorem{corollary}[theorem]{Corollary}
\newtheorem{lemma}[theorem]{Lemma}
\theoremstyle{definition}
\newtheorem{example}[theorem]{Example}
\newtheorem{c-example}[theorem]{Counter Example}

\newtheorem*{Beweis}{Proof}
\newtheorem{definition}[theorem]{Definition}
\newtheorem{punto}[theorem]{}
\newtheorem{remark}[theorem]{Remark}
\newtheorem{remarks}[theorem]{Remarks}
\CompileMatrices
\input{tcilatex}

\begin{document}

\title{On Linear Difference Equations over Rings and Modules \thanks{%
2000 MSC:\ Primary 16W30; Secondary 39A99.\newline
Key words and phrases: difference equations, linearly (bi)recursive
(bi)sequences, reversible sequences, Hopf algebras, dual coalgebra, rings,
modules.}}
\author{\textbf{Jawad Y. Abuhlail} \\
Mathematics Department\\
Birzeit University\\
P.O.Box 14, Birzeit - Palestine \\
jabuhlail@birzeit.edu}
\date{}
\maketitle

\begin{abstract}
In this note we develop a coalgebraic approach to the study of solutions of
linear difference equations over modules and rings. Some known results about
linearly recursive sequences over base fields are generalized to linearly
(bi)recursive (bi)sequences of modules over arbitrary commutative ground
rings.
\end{abstract}

\section*{Introduction}

Although the theory of linear difference equations over base fields is well
understood, the theory over arbitrary ground rings and modules is still
under development. It is becoming more interesting and is gaining
increasingly special importance mainly because of recent applications in
coding theory and cryptography (e.g. \cite{HN99}, \cite{KKMMN99}).

In a series of papers E. Taft et al. (e.g. \cite{PT80}, \cite{LT90}, \cite
{Taf95}) developed a coalgebraic aspect to the study of linearly recursive
sequences over fields. Moreover L. Gr\"{u}nenfelder et al. studied in (\cite
{GO93}, \cite{GK97}) the linearly recursive sequences over \emph{finite
dimensional }vector spaces. Linearly recursive (bi)sequences over arbitrary
rings and modules were studied intensively by A. Nechaev et al. (e.g. \cite
{Nec96}, \cite{KKMN95}, \cite{Nec93}), however the coalgebraic approach in
their work was limited to the field case. Generalization to the case of
arbitrary commutative ground rings was studied by several authors including
V. Kurakin (\cite{Kur94}, \cite{Kur00} $\&$ \cite{Kur2002}) and eventually
Abuhlail, Gomez-Torrecillas and Wisbauer \cite{AG-TW2000}.

In this note we develop a coalgebraic aspect to the study of solutions of
linear difference equations over \emph{arbitrary }rings and modules. For
some of our results we assume that the ground ring is artinian. Our results
generalize also previous results of us in \cite{AG-TW2000} and \cite[Kapitel
4]{Abu2001}. A standard reference for the theory of linearly recursive
sequences over rings and modules is the comprehensive work of A. Mikhalev et
al. \cite{KKMN95}. For the theory of Hopf algebras the reader may refer to
any of the classical references (e.g. \cite{Swe69}, \cite{Abe80} and \cite
{Mon93}).

With $R$ we denote a commutative ring with $1_{R}\neq 0_{R}$ and with $%
U(R)=\{r\in R|$ $r$ is invertible$\}$ the group of \emph{units} of $R.$ The
category of $R$-(bi)modules will be denoted by $\mathcal{M}_{R}.$ For an $R$%
-module $M,$ we call an $R$-submodule $K\subset M$ \emph{pure }(in the sense
of Cohn), if for every $R$-module $N$ the induced map $\iota _{k}\otimes 
\mathrm{id}_{N}:K\otimes _{R}N\rightarrow M\otimes _{R}N$ is injective.

For an $R$-algebra $A$ and an $A$-module $M,$ we call an $A$-submodule $%
K\subset M$ $R$-cofinite, if $M/K$ is f.g. in $\mathcal{M}_{R}.$ For an $R$%
-algebra $A$ we denote by $\mathcal{K}_{A}$ the class of $R$-cofinite
ideals. If $A$ is an $R$-algebra with $\mathcal{K}_{A}$ a filter, then we
define for every left $A$-module $M$ the\emph{\ finite dual} right $A$%
-module 
\begin{equation}
M^{\circ }:=\{f\in M^{\ast }\mid \text{ }\mathrm{Ke}(f)\supset IM\text{ for
some }A\text{-ideal }I\text{ with }A/I\text{ f.g.}\}.  \label{M0}
\end{equation}
With $\mathbb{N}$ resp. $\mathbb{Z}$ we denote the set of natural numbers
resp. the ring of integers. Moreover we set $\mathbb{N}_{0}:=\{0,1,2,3,...%
\}. $ For an $n\times n$ matrix $M$ over $R$ we denote the characteristic
polynomial with $\chi (M).$ The identity matrix of order $n$ over $R$ is
denoted by $E_{n}.$ For an $m\times n$ matrix $A$ and a $k\times l$ matrix $%
B,$ the \emph{Kronecker product} (\emph{tensor product}) of $A$ and $B$ is
the $mk\times nl$ matrix 
\begin{equation*}
A\otimes B:=\left[ 
\begin{array}{ccccc}
a_{11}\cdot B & a_{12}\cdot B & ... & ... & a_{1n}\cdot B \\ 
a_{21}\cdot B & a_{22}\cdot B & ... & ... & a_{2n}\cdot B \\ 
... & ... & ... & ... & ... \\ 
... & ... & ... & ... & ... \\ 
a_{m1}\cdot B & a_{m2}\cdot B & ... & ... & a_{mn}\cdot B
\end{array}
\right]
\end{equation*}

\section{Preliminaries}

Let $M$ be an $R$-module and 
\begin{equation}
M[\mathbf{x}]:=M[x_{1},...,x_{k}],\text{ }M[\mathbf{x},\mathbf{x}%
^{-1}]:=M[x_{1},x_{1}^{-1},...,x_{k},x_{k}^{-1}].  \label{Mx}
\end{equation}
We consider the \emph{polynomial ring} $R[\mathbf{x}]$ and the \emph{ring of
Laurent polynomials\ }$R[\mathbf{x},\mathbf{x}^{-1}]$ as commutative $R$%
-algebras with the usual multiplication and the usual unity. For every $R$%
-module $M,$ $M[\mathbf{x}]$ (resp. $M[\mathbf{x},\mathbf{x}^{-1}]$) is an $%
R[\mathbf{x}]$-module (resp. an $R[\mathbf{x},\mathbf{x}^{-1}]$-module) with
action induced from the $R$-module structure on $M$ and we have moreover
canonical $R$-module isomorphisms 
\begin{equation*}
M[\mathbf{x}]\simeq M\otimes _{R}R[\mathbf{x}]\simeq M^{(\mathbb{N}_{0}^{k})}%
\text{ and }M[\mathbf{x},\mathbf{x}^{-1}]\simeq M\otimes _{R}R[\mathbf{x},%
\mathbf{x}^{-1}]\simeq M^{(\mathbb{Z}^{k})}.
\end{equation*}
For $\mathbf{n}=(n_{1},...,n_{k})\in \mathbb{N}_{0}^{k}$ resp. $\mathbf{z}%
=(z_{1},...,z_{k})\in \mathbb{Z}^{k}$ we set $\mathbf{x}^{\mathbf{n}%
}:=x_{1}^{n_{1}}\cdot ...\cdot x_{k}^{n_{k}}$ resp. $\mathbf{x}^{\mathbf{z}%
}:=x_{1}^{z_{1}}\cdot ...\cdot x_{k}^{z_{k}}.$

\begin{punto}
\label{system}Let $M$ be an $R$-module, $\mathbf{l=}(l_{1},...,l_{k})\in 
\mathbb{N}_{0}^{k}$ and consider the \emph{system of linear difference
equations }(ab. SLDE) 
\begin{equation}
\begin{tabular}{lllll}
$x_{\mathbf{n}+(l_{1},0,...,0)}$ & $+$ & $\sum%
\limits_{i=1}^{l_{1}}p_{(1,l_{1}-i)}(\mathbf{n})x_{\mathbf{n}%
+(l_{1}-i,0,...,0)}$ & $=$ & $g_{1}(\mathbf{n}),$ \\ 
$x_{\mathbf{n}+(0,l_{2},0,...,0)}$ & $+$ & $\sum%
\limits_{i=1}^{l_{2}}p_{(2,l_{2}-i)}(\mathbf{n})x_{\mathbf{n}%
+(0,l_{2}-i,0,...,0)}$ & $=$ & $g_{2}(\mathbf{n}),$ \\ 
$...$ & $...$ & $...$ & $...$ & $...$ \\ 
$...$ & $...$ & $...$ & $...$ & $...$ \\ 
$x_{\mathbf{n}+(0,...,0,l_{k})}$ & $+$ & $\sum%
\limits_{i=1}^{l_{k}}p_{(k,l_{k}-i)}(\mathbf{n})x_{\mathbf{n}%
+(0,...,0,l_{k}-i)}$ & $=$ & $g_{k}(\mathbf{n}),$%
\end{tabular}
\label{SLDE}
\end{equation}
where the $p_{jl}$'s are $R$-valued functions and the $g_{j}$'s are $M$%
-valued functions defined for all $\mathbf{n}\in \mathbb{N}_{0}^{k}.$ If the 
$g_{j}$'s are identically zero, then (\ref{SLDE}) is said to be a \emph{%
homogenous }SLDE.\emph{\ }If the $p_{jl}$'s are constants, then (\ref{SLDE})
is said to be a SLDE \emph{with constant coefficients}.
\end{punto}

\begin{punto}
For an $R$-module $M$ and $k\geq 1$ let 
\begin{equation*}
\mathcal{S}_{M}^{<k>}:=\{u:\mathbb{N}_{0}^{k}\rightarrow M\}\simeq M^{%
\mathbb{N}_{0}^{k}}
\end{equation*}
be the $R$-module of $k$-\emph{sequences} over $M.$ If $M$ (resp. $k$) is
not mentioned, then we mean $M=R$ (resp. $k=1$). For $f(\mathbf{x}%
)=\sum\limits_{\mathbf{i}}a_{\mathbf{i}}\mathbf{x}^{\mathbf{i}}\in R[\mathbf{%
x}]$ and $w\in \mathcal{S}_{M}^{<k>}$ define 
\begin{equation}
f(\mathbf{x})\rightharpoonup w=u\in \mathcal{S}_{M}^{<k>},\text{ where }u(%
\mathbf{n}):=\sum_{\mathbf{i}}a_{\mathbf{i}}w(\mathbf{n+i})\text{ for all }%
\mathbf{n}\in \mathbb{N}_{0}^{k}\mathbb{.}  \label{f-th}
\end{equation}
With this action $\mathcal{S}_{M}^{<k>}$ is an $R[\mathbf{x}]$-module. For
subsets $I\subset R[\mathbf{x}]$ and $L\subset \mathcal{S}_{M}^{<k>}$
consider the annihilator submodules 
\begin{equation*}
\begin{tabular}{lll}
$\mathrm{An}_{\mathcal{S}_{M}^{<k>}}(I)$ & $=$ & $\{w\in \mathcal{S}%
_{M}^{<k>}\,|$ $f\rightharpoonup w=0$ for every $f\in I\},$ \\ 
$\mathrm{An}_{R[\mathbf{x}]}(L)$ & $=$ & $\{h\in R[\mathbf{x}%
]|\,h\rightharpoonup u=0$ for every $u\in L\}.$%
\end{tabular}
\end{equation*}
Note that $\mathrm{An}_{\mathcal{S}_{M}}^{<k>}(I)\subset \mathcal{S}%
_{M}^{<k>}$ is an $R[\mathbf{x}]$-submodule and $\mathrm{An}_{R[\mathbf{x}%
]}(L)\vartriangleleft R[\mathbf{x}]$ is an ideal.
\end{punto}

\begin{punto}
A polynomial $f(x)\in R[x]$ is called \emph{monic}, if its leading
coefficient is $1_{R}.$ For every monic polynomial $%
f(x)=x^{l}+a_{l-1}x^{l-1}+...+a_{1}x+a_{0}\in R[x],$ the \emph{companion
matrix}\textbf{\ }of $f$ is defined to be the $l\times l$ matrix 
\begin{equation}
S_{f}:=\left[ 
\begin{array}{ccccc}
0_{R} & 0_{R} & ... & 0_{R} & -a_{0} \\ 
1_{R} & 0_{R} & ... & 0_{R} & -a_{1} \\ 
0_{R} & 1_{R} & ... & 0_{R} & -a_{2} \\ 
... & ... & ... & ... & ... \\ 
0_{R} & 0_{R} & ... & 1_{R} & -a_{l-1}
\end{array}
\right]  \label{CM}
\end{equation}
$S_{f}$ is a matrix that has $f(x)$ as its characteristic polynomial as well
as its minimum polynomial (\cite[Theorem 4.18]{Jon73}).
\end{punto}

\begin{definition}
An ideal $I\vartriangleleft R[\mathbf{x}]$ will be called \emph{monic}, if
it contains a non-empty subset of monic polynomials 
\begin{equation}
\{f_{j}(x_{j})=x_{j}^{l_{j}}+a_{l_{j}-1}^{(j)}x_{j}^{l_{j}-1}+...+a_{1}^{(j)}x_{j}+a_{0}^{(j)}|%
\text{ }j=1,...,k\}.  \label{nP}
\end{equation}
In this case the polynomials (\ref{nP}) are called \emph{elementary
polynomials} and $(f_{1}(x_{1}),...,f_{k}(x_{k}))\vartriangleleft R[\mathbf{x%
}]$ an \emph{elementary ideal. }A monic polynomial $q(x)\in R[x]$ is called 
\emph{reversible}, if $q(0)\in U(R)$. An ideal $I\vartriangleleft R[\mathbf{x%
},\mathbf{x}^{-1}]$ will be called \emph{reversible}, if it contains a
subset of reversible polynomials $\{q_{1}(x_{1}),...,q_{k}(x_{k})\}.$
\end{definition}

\begin{punto}
Let $M$ be an $R$-module. We call $u\in \mathcal{S}_{M}^{<k>}$ a \emph{%
linearly recursive }$k$\emph{-sequence }(resp. a \emph{linearly birecursive }%
$k$\emph{-sequence}), if $\mathrm{An}_{R[\mathbf{x}]}(u)$ is a monic ideal
(resp. a reversible ideal). Note that a $k$-sequence $u\in \mathcal{S}%
_{M}^{<k>}$ is linearly recursive,\emph{\ }iff it's a solution of a
homogenous SLDE with constants coefficients of the form (\ref{SLDE}). If $%
\mathrm{An}_{R[\mathbf{x}]}(u)$ contains a set of monic polynomials $%
\{f_{1}(x_{1}),...,f_{k}(x_{k})\},$ where $f_{j}(x_{j})$ is of order $m_{j},$
$j=1,...,k,$ then these are called \emph{elementary characteristic
polynomials} of $u$ and $u$ is said to have \emph{order} $\mathbf{m}%
:=(m_{1},...,m_{k}).$ Characteristic polynomials of $u$ of least degree $%
n_{j},$ $j=1,...,k$ are called \emph{minimal polynomials }of $u$ and $%
\mathbf{n}:=(n_{1},...,n_{k})$ is called the \emph{rank} of $u.$ The subsets 
$\mathcal{L}_{M}^{<k>}\subseteq \mathcal{S}_{M}^{<k>}$ of linearly recursive 
$k$-sequences and $\mathcal{B}_{M}^{<k>}\subseteq \mathcal{S}_{M}^{<k>}$ of
linearly birecursive $k$-sequences are obviously $R[\mathbf{x}]$-submodules.
\end{punto}

\begin{punto}
(\cite[Page 170]{MN96})\label{lex} The \emph{lexicographical linear order} ($%
\preceq )$ on $\mathbb{N}_{0}^{k}$ is defined as follows: for $\mathbf{i}%
=(i_{1},...,i_{k})$ and $\mathbf{n}=(n_{1},...,n_{k})\in \mathbb{N}_{0}^{k}$
we say $\mathbf{i}\preceq \mathbf{n},$ if the first number in the sequence
of integers 
\begin{equation*}
(n_{1}+...+n_{k})-(i_{1}+...+i_{k}),\text{ }n_{1}-i_{1},...,n_{k}-i_{k}
\end{equation*}
that is different from zero is positive.

Let $M$ be an $R$-module, $\mathbf{F}:=\{f_{1}(x_{1}),...,f_{k}(x_{k})\}%
\subset R[\mathbf{x}]$ a subset of monic polynomials with $\mathrm{\deg }%
(f_{j}(x_{j}))=l_{j}$ for $j=1,...,k,$ $\mathbf{l}:=(l_{1},...,l_{k}),$ $%
\mathbf{1:}=(1,...,1),$ and $I_{\mathbf{F}}:=(f_{1},...,f_{k})%
\vartriangleleft R[\mathbf{x}].$ Note that the natural order ``$\leq "$ on $%
\mathbb{N}_{0}$ induces on $\mathbb{N}_{0}^{k}$ a \emph{partial order} and
we define the \emph{polyhedron }$\Pi _{\mathbf{F}}=\Pi (\mathbf{l}):=\{%
\mathbf{i}\in \mathbb{N}_{0}^{k}|$ $\mathbf{i}\leq \mathbf{l}-\mathbf{1}\}.$
The \emph{initial polyhedron of values}\textbf{\ }of $\omega \in \mathcal{S}%
_{M}^{<k>}$ is defined as $\omega (\Pi _{\mathbf{F}}):=\{\omega (\mathbf{i}%
)| $ $\mathbf{i}\in \Pi _{\mathbf{F}}\}.$ For $l=l_{1}\cdot ...\cdot l_{k}$
the points of the polyhedron $\Pi _{\mathbf{F}}$ build a chain $\mathbf{0}=%
\mathbf{i}_{0}\preceq \mathbf{i}_{1}\preceq ...\preceq \mathbf{i}_{l-1}$ and
we can write $\omega (\Pi _{\mathbf{F}})$ as an \emph{initial vector\ of
values}\textbf{\ }$(\omega (\mathbf{0}),\omega (\mathbf{i}_{1}),...,\omega (%
\mathbf{i}_{l-1}))\in M^{l}.$

Let $\omega \in \mathrm{An}_{\mathcal{S}%
_{M}^{<k>}}(f_{1}(x_{1}),...,f_{k}(x_{k})),$ where $f_{j}(x_{j})$ is monic
for $j=1,...,n$ and write for every $\mathbf{n}=(n_{1},...,n_{k})\in \mathbb{%
N}_{0}^{k}:$%
\begin{equation*}
x_{j}^{n_{j}}=h_{j}(x_{j})f_{j}(x_{j})+r_{j}(x_{j}),\text{ where }\mathrm{%
\deg }(r_{j}(x_{j}))<l_{j}.
\end{equation*}
If we set 
\begin{equation*}
g^{(\mathbf{n})}(\mathbf{x}{\normalsize ):=}\prod_{j=1}^{k}r_{j}(x_{j})=\sum%
\limits_{\mathbf{i}\in \Pi _{\mathbf{F}}}a_{\mathbf{i}}^{(\mathbf{n}%
{\normalsize )}}\mathbf{x}{\normalsize ^{\mathbf{i}}}\text{ and }v:=\mathbf{x%
}^{\mathbf{n}}\rightharpoonup \omega =g^{(\mathbf{n})}(\mathbf{x}%
)\rightharpoonup \omega ,
\end{equation*}
then 
\begin{equation*}
\omega (\mathbf{n}{\normalsize )=v(\mathbf{0})=\sum\limits_{\mathbf{i}\in
\Pi _{\mathbf{F}}}a_{\mathbf{i}}^{(\mathbf{n})}\omega (}\mathbf{i}%
{\normalsize )}\text{ for every }\mathbf{n}\in \mathbb{N}_{0}^{k}%
{\normalsize .}
\end{equation*}
Consequently $\omega $ is completely determined by the initial polyhedron of
values $\omega (\Pi _{\mathbf{F}}).$ For $\mathbf{t}\in \Pi _{\mathbf{F}}$
define the sequence $e_{\mathbf{t}}^{\mathbf{F}}\in \mathrm{An}_{\mathcal{S}%
_{R}^{<k>}}(I_{\mathbf{F}})$ with initial polyhedron of values $e_{\mathbf{t}%
}^{\mathbf{F}}(\mathbf{i})=\delta _{\mathbf{i},\mathbf{t}}$ for all $\mathbf{%
i}\in \Pi _{\mathbf{F}}.$ The sequence $e_{\mathbf{l}-\mathbf{1}}^{\mathbf{F}%
}$ is called the \emph{impulse sequence\ }of $\mathrm{An}_{\mathcal{S}%
_{R}^{<k>}}(I_{\mathbf{F}}).$
\end{punto}

\section*{Examples}

We give now some examples of linearly recursive sequences. For more examples
the reader may refer to \cite{KKMN95}.

\begin{example}
(\emph{Geometric progression}). Let $M$ be an $R$-module, $m\in M,$ $r\in R$
and consider $w\in \mathcal{S}_{M}$ given by 
\begin{equation*}
w(n{\normalsize ):=r^{n}m}\text{ for every }{\normalsize n\in \mathbb{N}_{0}.%
}
\end{equation*}
Then $w\in \mathcal{L}_{M}$ with initial condition $w(0)=m$ and elementary
characteristic polynomial $f(x)=x-r.$ Moreover $\mathrm{An}%
_{R[x]}(w)=R[x](x-r)+R[x]\mathrm{An}_{R}(r).$
\end{example}

\begin{example}
(\emph{Arithmetic progression}). Let $M$ be an $R$-module, $\{p,q\}\subset M$
and consider $w\in \mathcal{S}_{M}$ given by 
\begin{equation*}
w(n{\normalsize ):=p+nq}\text{ for every }n\in {\normalsize \mathbb{N}_{0}.}
\end{equation*}
Then $w\in \mathcal{L}_{M}$ with initial vector $(p,p+q)$ and elementary
characteristic polynomial $f(x)=(x-1)^{2}.$ If $\mathrm{An}_{R}(q)=0,$ then $%
f(x)$ is a unique minimal polynomial of $w.$ If $r\in \mathrm{An}_{R}(q),$
then $f_{r}(x)=(x-1)^{2}+r(x-1)$ is another minimal polynomial of $w.$
\end{example}

\begin{remark}
An example of a \emph{non} linearly recursive sequences over $\mathbb{Z}$ is
the sequence of prime positive numbers $\{2,3,5,7,...\}.$
\end{remark}

\begin{example}
Let $E=\{f_{1}(x),...,f_{k}(x)\}\subset R[x]$ be a subset of monic
polynomials.

\begin{enumerate}
\item  Let $M$ be an $R$-module, $u_{i}\in \mathrm{An}_{\mathcal{S}%
_{M}}(f_{i})$ for $i=1,...,k$ and consider $u:=u_{1}\overset{\cdot }{+}...%
\overset{\cdot }{+}u_{k}\in \mathcal{S}_{M}^{<k>}$ defined by $u(\mathbf{n}%
)=u_{1}(n_{1})+...+u_{k}(n_{k}).$ Then $u\in \mathrm{An}_{\mathcal{S}%
_{M}^{<k>}}(g_{1}(x_{1}),...,g_{k}(x_{k})),$ where for $i=1,...,k:$%
\begin{equation}
g_{i}(x_{i})= 
\begin{cases}
f_{i}(x_{i}), & f_{i}(1_{R})=0_{R} \\ 
f_{i}(x_{i})(x_{i}-1_{R}), & \text{otherwise.}
\end{cases}
\label{gi}
\end{equation}

\item  Let $M_{1},...,M_{k}$ be $R$-modules, $u_{i}\in \mathrm{An}_{\mathcal{%
S}_{M_{i}}}(f_{i})$ for $i=1,...,k,$ $M:=M_{1}\oplus ...\oplus M_{k}$ and
consider $u\in \mathcal{S}_{M}^{<k>}$ defined by $u(\mathbf{n}%
):=(u_{1}(n_{1}),...,u_{k}(n_{k})).$ Then $u\in \mathrm{An}_{\mathcal{S}%
_{M}^{<k>}}(g_{1}(x_{1}),...,g_{k}(x_{k})),$ where the $g_{i}$'s are defined
as in (\ref{gi}).

\item  Let $u_{i}\in \mathrm{An}_{\mathcal{S}_{R}}(f_{i})$ for $i=1,...,k$
and consider $u\in \mathcal{S}_{R}^{<k>}$ defined by $u(\mathbf{n}%
):=u_{1}(n_{1})\cdot ...\cdot u_{k}(n_{k}).$ Then $u\in \mathrm{An}_{%
\mathcal{S}_{R}^{<k>}}(f_{_{1}}(x_{1}),...,f_{k}(x_{k}))$ and 
\begin{equation*}
\mathrm{An}_{\mathcal{S}_{R}^{<k>}}(f_{1}(x_{1}),...,f_{k}(x_{k}))\simeq 
\mathrm{An}_{\mathcal{S}_{R}}(f_{1})\otimes _{R}...\otimes _{R}\mathrm{An}_{%
\mathcal{S}_{R}}(f_{k}).
\end{equation*}

\item  Let $M_{1},...,M_{k}$ be $R$-modules, $u_{i}\in \mathrm{An}_{\mathcal{%
S}_{M_{i}}}(f_{i})$ for $i=1,...,k,$ $M:=M_{1}\otimes _{R}...\otimes
_{R}M_{k}$ and consider $u\in \mathcal{S}_{M}^{<k>}$ defined by $u(\mathbf{n}%
):=u_{1}(n_{1})\otimes ...\otimes u_{k}(n_{k}).$ Then $u\in \mathrm{An}_{%
\mathcal{S}_{M}^{<k>}}(f_{1}(x_{1}),...,f_{k}(x_{k}))$ and 
\begin{equation*}
\mathrm{An}_{\mathcal{S}_{M}^{<k>}}(f_{1}(x_{1}),...,f_{k}(x_{k}))\simeq 
\mathrm{An}_{\mathcal{S}_{M_{1}}}(f_{1})\otimes _{R}...\otimes _{R}\mathrm{An%
}_{\mathcal{S}_{M_{k}}}(f_{k}).
\end{equation*}
\end{enumerate}
\end{example}

\section*{Admissible $R$-bialgebras and Hopf $R$-algebras}

For every $R$-coalgebra $(C,\Delta _{C},\varepsilon _{C})$ there is a \emph{%
dual }$R$\emph{-algebra} $C^{\ast }:=\mathrm{Hom}_{R}(C,R)$ with
multiplication the so called \emph{convolution product} 
\begin{equation*}
(f\star g)(c):=\sum f(c_{1})g(c_{2})\text{ for all }f,g\in C^{\ast },\text{ }%
c\in C
\end{equation*}
and unity $\varepsilon _{C}.$ Although every algebra $A$ has a \emph{dual
coalgebra, }if the ground ring is hereditary noetherian (e.g. a field), the
existence of dual coalgebras\emph{\ }of algebras over an arbitrary
commutative ground rings is not guaranteed!! One way to handle this problem
is to restrict the class of $R$-algebras, for which the dual $R$-coalgebras
are defined.

\begin{definition}
Let $A$ be an $R$-algebra (resp. an $R$-bialgebra, a Hopf $R$-algebra). Then
we call $A:$

\begin{enumerate}
\item  an $\alpha $\emph{-algebra }(resp. an $\alpha $\emph{-bialgebra,} a 
\emph{Hopf }$\alpha $\emph{-algebra}),\emph{\ }if $\mathcal{K}_{A}$ is a
filter and $A^{\circ }\subset R^{A}$ is pure.

\item  \emph{cofinitary,}\textbf{\ }if $\mathcal{K}_{A}$ is a \emph{filter}
and for every $I\in \mathcal{K}_{A}$ there exists an $A$-ideal $\overline{I}%
\subseteq I$ with $A/\overline{I}$ f.g. and projective.
\end{enumerate}
\end{definition}

\begin{punto}
\label{zul-filt}{\normalsize \ }Let $H$ be an $R$-bialgebra and consider the
class of $R$-cofinite $H$-ideals $\mathcal{K}_{H}.$ We call $H$ an \emph{%
admissible }$R$\emph{-bialgebra}$,$ if $H$ is cofinitary and $\mathcal{K}%
_{H} $ satisfies the following axioms: 
\begin{equation}
\begin{tabular}{ll}
(A1) & $\forall $ $I,J\in \mathcal{K}_{H}$ there exists $L\in \mathcal{K}%
_{H},$ s.t. $\Delta _{H}(L)\subseteq I\otimes _{R}H+H\otimes _{R}J$%
\end{tabular}
\label{A1}
\end{equation}
and 
\begin{equation}
\begin{tabular}{ll}
(A2) & $\exists $ $I\in \mathcal{K}_{H},$ s.t. $\mathrm{Ke}(\varepsilon
_{H})\supset I.$%
\end{tabular}
\label{A2}
\end{equation}
We call a Hopf $R$-algebra $H$ an \emph{admissible Hopf }$R$\emph{-algebra},
if $H$ is cofinitary, $\mathcal{K}_{H}\ $satisfies (A1), (A2) and 
\begin{equation}
\begin{tabular}{ll}
(A3) & $\text{for every }I\in \mathcal{K}_{H}\text{ there exists }J\in 
\mathcal{K}_{H},\text{ s.t. }S_{H}(J)\subseteq I.$%
\end{tabular}
\label{A3}
\end{equation}
\end{punto}

\begin{remark}
\label{ad}It follows from the proof of \cite[Proposition 4.2.]{AG-TL2001},
that every cofinitary $R$-algebra (resp. $R$-bialgebra, Hopf $R$-algebra) is
an $\alpha $-algebra (resp. an $\alpha $-bialgebra, a Hopf $\alpha $%
-algebra). By (\cite[Lemma 2.5.6.]{Abu2001}) every cofinitary bialgebra
(Hopf algebra) over a \emph{noetherian }ground ring is admissible.
\end{remark}

\begin{proposition}
\label{cofinitary}\emph{(\cite[Proposition 2.4.13, Proposition 2..5.7]
{Abu2001})}

\begin{enumerate}
\item  If $A$ is a cofinitary $R$-algebra, then $A^{\circ }$ is an $R$%
-coalgebra. If $H$ is an admissible $R$-bialgebra \emph{(}resp. an
admissible Hopf $R$-algebra\emph{)}, then $H^{\circ }$ is an $R$-bialgebra 
\emph{(}resp. a Hopf $R$-algebra\emph{)}.

\item  Let $R$ be noetherian. If $A$ is an $\alpha $-algebra \emph{(}resp.
an $\alpha $-bialgebra, a Hopf $\alpha $-algebra\emph{)}, then $A^{\circ }$
is an $R$-coalgebra \emph{(}resp. an $R$-bialgebra, a Hopf $R$-algebra\emph{)%
}.
\end{enumerate}
\end{proposition}

\begin{proposition}
\label{sg}Let $A$ be an $\alpha $-algebra \emph{(}resp. an $\alpha $%
-bialgebra, a Hopf $\alpha $-algebra\emph{),} $B$ a cofinitary $R$-algebra 
\emph{(}resp. $R$-bialgebra, Hopf $R$-algebra\emph{)} and consider the
canonical map $\sigma :A^{\circ }\otimes _{R}B^{\circ }\rightarrow (A\otimes
_{R}B)^{\circ }.$ Then:

\begin{enumerate}
\item  $\sigma $ is injective.

\item  If $R$ is noetherian, then $\sigma $ is an isomorphism of $R$%
-coalgebras \emph{(}resp. $R$-bialgebras, Hopf $R$-algebras\emph{)}.
\end{enumerate}
\end{proposition}

\begin{Beweis}
\begin{enumerate}
\item  The proof is along the lines of the proof of \cite[Proposition 5]
{Kur2002}.

\item  The proof is along the lines of the proof of \cite[Theorem 4.10]
{AG-TL2001}.
\end{enumerate}
\end{Beweis}

The proof of \cite[Lemma 4.12]{AG-TL2001}\emph{\ }can be generalized to get

\begin{lemma}
\label{q-q}For any set of reversible polynomials $%
\{q_{1}(x_{1}),...,q_{k}(x_{k})\}\subseteq R[\mathbf{x}]$ we have an
isomorphism of $R$-algebras 
\begin{equation*}
R[\mathbf{x}]/(q_{1}(x_{1}),...,q_{k}(x_{k}))\simeq R[\mathbf{x},\mathbf{x}%
^{-1}]/(q_{1}(x_{1}),...,q_{k}(x_{k})).
\end{equation*}
\end{lemma}

\begin{lemma}
\label{cof}\emph{(\cite[Proposition 1]{Kur2002})} Let $R$ be an \emph{%
arbitrary} commutative ring.

\begin{enumerate}
\item  An ideal $I\vartriangleleft R[\mathbf{x}]$ is $R$-cofinite, iff it's
monic. Consequently every $R$-cofinite $R[\mathbf{x}]$-ideal contains an
ideal $\overline{I}\vartriangleleft R[\mathbf{x}],$ such that $R[\mathbf{x}]/%
\overline{I}$ is free of finite rank. In particular $R[\mathbf{x}]$ is
cofinitary.

\item  An ideal $I\vartriangleleft R[\mathbf{x},\mathbf{x}^{-1}]$ is $R$%
-cofinite, iff it's reversible. Consequently every $R$-cofinite $R[\mathbf{x}%
,\mathbf{x}^{-1}]$-ideal contains an ideal $\overline{I}\vartriangleleft R[%
\mathbf{x},\mathbf{x}^{-1}],$ such that $R[\mathbf{x},\mathbf{x}^{-1}]/%
\overline{I}$ is free of finite rank. In particular $R[\mathbf{x},\mathbf{x}%
^{-1}]$ is cofinitary.
\end{enumerate}
\end{lemma}

\section{Linearly (bi)recursive sequences}

In this section we study the \emph{linearly }(\emph{bi})\emph{recursive }$k$%
\emph{-sequences} over $R$-modules, where $R$ is an arbitrary commutative
ground ring.

\begin{punto}
\label{kokomm}Let $(G,\mu _{G},e_{G})$ be a (commutative) monoid.
Considering the elements of the basis $G$ as \emph{group-like elements,} the
monoid algebra $RG$ becomes a (commutative) cocommutative $R$-bialgebra $%
(RG,\mu ,\eta ,\Delta _{g},\varepsilon _{g}),$ where 
\begin{equation*}
\Delta _{g}(x)=x\otimes x\text{ and }\varepsilon _{g}(x)=1_{R}\text{ for
every }x\in G.
\end{equation*}
If $G$ is a group, then $RG$ is a Hopf $R$-algebra with antipode 
\begin{equation*}
S_{g}:RG\rightarrow RG,\text{ }x\mapsto x^{-1}\text{ for every }x\in G.
\end{equation*}
\end{punto}

\begin{punto}
\textbf{Bialgebra structures on }$R[\mathbf{x}]$\label{rx-co}. Consider the
commutative monoid $G$ generated by $\{x_{j}\mid j=1,...,k\}.$ Then $R[%
\mathbf{x}]=RG$ has the structure of a \emph{commutative cocommutative}%
{\normalsize \ }$R$-bialgebra $R[\mathbf{x};g]=(R[\mathbf{x}],\mu ,\eta
,\Delta _{g},\varepsilon _{g}),$ where $\mu $ is the usual multiplication, $%
\eta $ is the usual unity and 
\begin{equation*}
\begin{tabular}{llllllll}
$\Delta _{g}:$ & $R[\mathbf{x}]$ & $\rightarrow $ & $R[\mathbf{x}]\otimes
_{R}R[\mathbf{x}],$ & $x_{j}^{n}$ & $\mapsto $ & $x_{j}^{{n}}\otimes x_{j}^{{%
n}},$ & $\forall $ $n\geq 0,$ $j=1,...,k,$ \\ 
$\varepsilon _{g}:$ & $R[\mathbf{x}]$ & $\rightarrow $ & $R,$ & $x_{j}^{n}$
& $\mapsto $ & $1_{R},$ & $\forall $ $n\geq 0,$ $j=1,...,k.$%
\end{tabular}
\end{equation*}
On the other hand $R[\mathbf{x};p]=(R[\mathbf{x}],\mu ,\eta ,\Delta
_{g},\varepsilon _{g})$ is a \emph{commutative cocommutative\ }Hopf $R$%
-algebra, where $\mu $ is the usual multiplication, $\eta $ is the usual
unity and 
\begin{equation*}
\begin{tabular}{lllllllll}
$\Delta _{p}$ & $:$ & $R[\mathbf{x}]$ & $\rightarrow $ & $R[\mathbf{x}%
]\otimes _{R}R[\mathbf{x}],$ & $x_{j}^{{n}}$ & $\mapsto $ & $%
\sum\limits_{t=0}^{{n}}\binom{n}{t}$ $x_{j}^{{t}}\otimes x_{j}^{n{-t}},$ & $%
\forall $ $n\geq 0,$ $j=1,...,k,$ \\ 
$\varepsilon _{p}$ & $:$ & $R[\mathbf{x}]$ & $\rightarrow $ & $R,$ & $x_{j}^{%
{n}}$ & $\mapsto $ & $\delta _{n,0},$ & $\forall $ $n\geq 0,$ $j=1,...,k,$
\\ 
$S_{p}$ & $:$ & $R[\mathbf{x}]$ & $\rightarrow $ & $R[\mathbf{x}],$ & $%
x_{j}^{{n}}$ & $\mapsto $ & $(-1)^{n}x_{j}^{n},$ & $\forall $ $n\geq 0,$ $%
j=1,...,k.$%
\end{tabular}
\end{equation*}
\end{punto}

\begin{remarks}
\begin{enumerate}
\item  Let $R$ be an integral domain, then it follows by \cite[Theorem
1.3.6.]{Gru69} that for every set $G,$ the class of group-like elements of
the $R$-coalgebra $RG$ is $G$ itself. Then one can show as in the field case 
\cite{CG93}, that $R[\mathbf{x};g]$ and $R[\mathbf{x};p]$ are the only
possible $R$-bialgebra structures on $R[\mathbf{x}]$ with the usual
multiplication and the usual unity.

\item  The $R$-bialgebra $R[\mathbf{x};g]$ has no antipode, because the
group-like elements in a Hopf $R$-algebra should be invertible.
\end{enumerate}
\end{remarks}

The proof of the following result depends mainly on arguments of 
\cite[Theorem 2]{Kur2002}:

\begin{proposition}
\label{rx-admiss}Let $R$ be an \emph{arbitrary} commutative ring. Then $R[%
\mathbf{x};g]$ is an admissible $R$-bialgebra and $R[\mathbf{x};p]$ is an
admissible Hopf $R$-algebra. Hence $R[\mathbf{x};g]^{\circ }$ is an $R$%
-bialgebra and $R[\mathbf{x};p]^{\circ }$ is a Hopf $R$-algebra.
\end{proposition}

\begin{Beweis}
Denote with $(R[\mathbf{x}],\Delta ,\varepsilon )$ either of the cofinitary $%
R$-bialgebras $R[\mathbf{x};g]$ and $R[\mathbf{x};p].$ Let $%
I,J\vartriangleleft R[\mathbf{x}]$ be $R$-cofinite ideals and assume
w.l.o.g. that $R[\mathbf{x}]/I$ and $R[\mathbf{x}]/J$ are free of finite
rank (see Lemma \ref{cof}). Let $\beta $ be a basis of the free $R$-module $%
B:=R[\mathbf{x}]/I\otimes _{R}R[\mathbf{x}]/J$ and consider the $R$-algebra
morphism $\overline{\Delta }:=(\pi _{I}\otimes \pi _{J})\circ \Delta :R[%
\mathbf{x}]\rightarrow R[\mathbf{x}]/I\otimes _{R}R[\mathbf{x}]/J.$ For $%
j=1,...,k$ let $M_{j}$ be the matrix of the $R$-linear map 
\begin{equation*}
T_{j}:B\rightarrow B,\text{ }b\mapsto \overline{\Delta }(x_{j})b
\end{equation*}
w.r.t. $\beta $ and $\chi _{j}(\lambda )$ its characteristic polynomial.
Then $\chi _{j}(\overline{\Delta }(x_{j}))=0$ for $j=1,...,k.$ Since $%
\overline{\Delta }$ is an $R$-algebra morphism, it follows that $\chi
_{j}(x_{j})\in \mathrm{Ke}(\overline{\Delta })=\Delta ^{-1}(I\otimes _{R}R[%
\mathbf{x}]+R[\mathbf{x}]\otimes _{R}J)$ for $j=1,...,k.$ If we set $%
L:=(\chi _{1}(x_{1}),...,\chi _{k}(x_{k}))\vartriangleleft R[\mathbf{x}],$
then $\Delta (L)\subseteq I\otimes _{R}R[\mathbf{x}]+R[\mathbf{x}]\otimes
_{R}J,$ i.e. $\mathcal{K}_{R[\mathbf{x}]}$ satisfies axiom (\ref{A1}). Note
that $R[\mathbf{x}]/\mathrm{Ke}(\varepsilon )\simeq R,$ hence $\mathcal{K}%
_{R[\mathbf{x}]}$ satisfies axiom (\ref{A2}). Consequently $R[\mathbf{x};g]$
and $R[\mathbf{x};p]$ are admissible $R$-bialgebras. Consider now the Hopf $%
R $-algebra $R[\mathbf{x};p]$ with the \emph{bijective }antipode $S_{p}.$
For every ideal $I\vartriangleleft R[\mathbf{x}],$ $S_{p}^{-1}(I)%
\vartriangleleft R[\mathbf{x};p]$ is an ideal and we have an isomorphism of $%
R$-modules $R[\mathbf{x}]/S_{p}^{-1}(I)\simeq R[\mathbf{x}]/I,$ hence $%
\mathcal{K}_{R[\mathbf{x};p]}$ satisfies axiom (\ref{A3}). Consequently $R[%
\mathbf{x};p]$ is an admissible Hopf $R$-algebra. The last statement follows
now by Proposition \ref{cofinitary}.$\blacksquare $
\end{Beweis}

If $M$ is an arbitrary $R$-module, then we have obviously an isomorphism of $%
R[\mathbf{x}]$-modules 
\begin{equation}
\Phi _{M}:M[\mathbf{x}]^{\ast }\rightarrow \mathcal{S}_{M^{\ast }}^{<k>},%
\text{ }\varkappa \mapsto \lbrack \mathbf{n}\mapsto \lbrack m\mapsto
\varkappa (m\mathbf{x}^{\mathbf{n}})]]  \label{folge-M}
\end{equation}
with inverse $u\mapsto \lbrack m\mathbf{x}^{\mathbf{n}}\mapsto u(\mathbf{n)(}%
m)].$

\begin{proposition}
\label{mxn-d-l}Let $M$ be an $R$-module. Then \emph{(\ref{folge-M}) }induces
an isomorphism of $R[\mathbf{x}]$-modules 
\begin{equation}
M[\mathbf{x}]^{\circ }\simeq {\normalsize \mathcal{L}}_{M^{\ast }}^{<k>}.
\label{Mx0=L_Mk}
\end{equation}
\end{proposition}

\begin{Beweis}
Consider the $R[\mathbf{x}]$-module isomorphism $M[\mathbf{x}]^{\ast }%
\overset{\Phi _{M}}{\simeq }\mathcal{S}_{M^{\ast }}^{<k>}$ (\ref{folge-M}).
Let $\varkappa \in M[\mathbf{x}]^{\circ }.$ Then there exists an $R$%
-cofinite $R[\mathbf{x}]$-ideal $I,$ such that $I\rightharpoonup \varkappa
=0.$ So $I\rightharpoonup \Phi (\varkappa )=\Phi (I\rightharpoonup \varkappa
)=0,$ i.e. $I\subset \mathrm{An}_{R[\mathbf{x}]}(\Phi (\varkappa )).$ By
Lemma \ref{cof} (1) $I$ is monic, i.e. $\Phi (\varkappa )\in \mathcal{L}%
_{M^{\ast }}^{<k>}.$

On the other hand, let $u\in \mathcal{L}_{M^{\ast }}^{<k>}.$ By definition $%
J:=\mathrm{An}_{R[\mathbf{x}]}(u)$ is a monic ideal and it follows by Lemma 
\ref{cof} (1) that $J\vartriangleleft R[\mathbf{x}]$ is $R$-cofinite. For $%
\varkappa :=\Phi ^{-1}(u)$ we have $J\rightharpoonup \varkappa
=J\rightharpoonup \Phi ^{-1}(u)=\Phi ^{-1}(J\rightharpoonup u)=0,$ i.e. $%
\varkappa \in M[\mathbf{x}]^{\circ }.\blacksquare $
\end{Beweis}

\begin{punto}
\textbf{The coalgebra structure on }$\mathcal{L}^{<k>}.$\label{LR-co}

By Lemma \ref{cof} (1) $(R[\mathbf{x}],\mu ,\eta )$ is a cofinitary $R$%
-algebra, where $\mu $ is the usual multiplication and $\eta $ is the usual
unity. Hence $(R[\mathbf{x}]^{\circ },\mu ^{\circ },\eta ^{\circ })$ is (by
Proposition \ref{cofinitary}) an $R$-coalgebra, where 
\begin{equation*}
\begin{tabular}{lllllll}
$\mu ^{\circ }:$ & $R[\mathbf{x}]^{\circ }\rightarrow $ & $R[\mathbf{x}%
]^{\circ }\otimes _{R}R[\mathbf{x}]^{\circ },$ & $f$ & $\mapsto $ & $%
[x_{i}^{s}\otimes x_{j}^{t}\mapsto f(x_{i}^{s}x_{j}^{t}),$ $s,t\geq
0,i,j=1,...,k],$ &  \\ 
$\eta ^{\circ }:$ & $R[\mathbf{x}]^{\circ }\rightarrow $ & $R,$ & $f$ & $%
\mapsto $ & $\text{ }f(1_{R}).$ & 
\end{tabular}
\end{equation*}
So $\mathcal{L}^{<k>}\simeq R[\mathbf{x}]^{\circ }$ has the structure of an $%
R$-coalgebra with counity 
\begin{equation}
\varepsilon _{\mathcal{L}^{<k>}}:\mathcal{L}^{<k>}\rightarrow R,\text{ }%
u\mapsto u(\mathbf{0}).  \label{eps-L}
\end{equation}
and comultiplication described as follows (see \cite[Proposition 14.16]
{KKMN95}):

Let $u\in \mathcal{L}^{<k>},$ $\{f_{1}(x_{1}),...,f_{k}(x_{k})\}\subseteq 
\mathrm{An}_{R[\mathbf{x}]}(u)$ a subset of elementary characteristic
polynomials with $\mathrm{\deg }(f_{j}(x_{j}))=l_{j}$ and $\mathbf{l:}%
=(l_{1},...,l_{k}).$ So we have for all $\mathbf{n},\mathbf{i}\in \mathbb{N}%
_{0}^{k}:$%
\begin{equation*}
u(\mathbf{n}+\mathbf{i})=(\mathbf{x}^{\mathbf{i}}\rightharpoonup u)(\mathbf{n%
})=(\sum\limits_{\mathbf{t}\leq \mathbf{l}-\mathbf{1}}(\mathbf{x}^{\mathbf{i}%
}\rightharpoonup u)(\mathbf{t})\cdot e_{\mathbf{t}}^{\mathbf{F}})(\mathbf{n}%
)=\sum\limits_{\mathbf{t}\leq \mathbf{l}-\mathbf{1}}(\mathbf{x}^{\mathbf{t}%
}\rightharpoonup u)(\mathbf{i})\cdot e_{\mathbf{t}}^{\mathbf{F}}(\mathbf{n}).
\end{equation*}
The comultiplication of $\mathcal{L}^{<k>}$ is given then by 
\begin{equation}
\Delta _{\mathcal{L}^{<k>}}:\mathcal{L}^{<k>}\rightarrow \mathcal{L}^{<k>}{%
\otimes _{R}}\text{ }\mathcal{L}^{<k>},\text{ }u\mapsto \sum\limits_{\mathbf{%
t}\leq \mathbf{l}-\mathbf{1}}(\mathbf{x}^{\mathbf{t}}\rightharpoonup
u)\otimes e_{\mathbf{t}}^{\mathbf{F}}.  \label{delta-L}
\end{equation}
\end{punto}

\begin{example}
Consider the \emph{Fibonacci sequence} $\digamma =(0,1,1,2,3,5,...).$
Clearly $\digamma $ is given by 
\begin{equation*}
\digamma (0)=0,\text{ }\digamma (1)=1,\text{ }\digamma (n+2)=\digamma
(n+1)+\digamma (n)\text{ for all }n\geq 0,
\end{equation*}
i.e. $\digamma \in \mathcal{L}_{\mathbb{Z}}$ with initial vector $(0,1)$ and
elementary characteristic polynomial $f(x)=x^{2}-x-1\in \mathbb{Z}[x].$ By (%
\ref{delta-L}) one can easily calculate 
\begin{equation*}
\Delta _{\mathcal{L}_{\mathbb{Z}}}(\digamma )=\digamma \otimes _{\mathbb{Z}%
}(x\rightharpoonup \digamma )+(x\rightharpoonup \digamma )\otimes _{\mathbb{Z%
}}\digamma -\digamma \otimes _{\mathbb{Z}}\digamma .
\end{equation*}
\end{example}

\begin{punto}
\textbf{The }$R$\textbf{-bialgebra }$(\mathcal{L}_{R}^{<k>};g).$\label{L-Had}
Consider the $R$-bialgebra $R[\mathbf{x};g].$ Then $\mathcal{S}^{<k>}\simeq
R^{\mathbb{N}_{0}^{k}}\simeq R[\mathbf{x};g]^{\ast }$ is an $R$-algebra with
multiplication given by the \emph{Hadamard product} 
\begin{equation}
\ast _{g}:\mathcal{S}^{<k>}\otimes _{R}\mathcal{S}^{<k>}\rightarrow \mathcal{%
S}^{<k>},\text{ }u\otimes v\mapsto \lbrack \mathbf{n}\mapsto u(\mathbf{n})v(%
\mathbf{n})]  \label{Hadamard}
\end{equation}
and the unity 
\begin{equation}
\eta _{g}:R\rightarrow \mathcal{S}{\normalsize ^{<k>}},\text{ }1_{R}\mapsto
\lbrack \mathbf{n}\mapsto 1_{R}]\text{ for every }\mathbf{n}\in \mathbb{N}%
_{0}^{k}.  \label{eta-g}
\end{equation}
By Propositions \ref{rx-admiss} and \ref{mxn-d-l} $(\mathcal{L}%
_{R}^{<k>};g)\simeq R[\mathbf{x};g]^{\circ }$ has the structure of an $R$%
-bialgebra with the coalgebra structure described in \ref{LR-co}, the
Hadamard product (\ref{Hadamard}) and the unity (\ref{eta-g}).
\end{punto}

\begin{punto}
\textbf{The Hopf }$R$\textbf{-algebra }$(\mathcal{L}_{R}^{<k>};p).$\label%
{L-Hur} Consider the Hopf $R$-algebra $R[\mathbf{x};p].$ Then $\mathcal{S}%
^{<k>}\simeq R^{\mathbb{N}_{0}^{k}}\simeq R[\mathbf{x};p]^{\ast }$ is an $R$%
-algebra with multiplication given by the \emph{Hurwitz product} 
\begin{equation}
\ast _{p}:\mathcal{S}^{<k>}\otimes _{R}\mathcal{S}^{<k>}\rightarrow \mathcal{%
S}^{<k>},\text{ }u\otimes v\mapsto \lbrack \mathbf{n}\mapsto \sum_{\mathbf{t}%
\leq \mathbf{n}}\binom{\mathbf{n}}{\mathbf{t}}u(\mathbf{t})v(\mathbf{n}-%
\mathbf{t})]  \label{Hurwiz}
\end{equation}
and the unity 
\begin{equation}
\eta _{p}:R\rightarrow \mathcal{S}{\normalsize ^{<k>},}\text{ }1_{R}\mapsto
\lbrack \mathbf{n}\mapsto \delta _{\mathbf{n},\mathbf{0}}]\text{ for every }%
\mathbf{n}\in \mathbb{N}_{0}^{k}.  \label{eta-p}
\end{equation}
By Propositions \ref{rx-admiss} and \ref{mxn-d-l} $(\mathcal{L}%
_{R}^{<k>};p)\simeq R[\mathbf{x};p]^{\circ }$ has the structure of a Hopf $R$%
-algebra with the coalgebra structure described in \ref{LR-co}, the Hurwitz
product (\ref{Hurwiz}), the unity (\ref{eta-p}) and the antipode 
\begin{equation*}
S_{\mathcal{L}^{<k>}}:\mathcal{L}^{<k>}\rightarrow \mathcal{L}^{<k>},\text{ }%
u\mapsto \lbrack {\normalsize \mathbf{i}\mapsto (-1)^{\mathbf{i}}u(\mathbf{i}%
)]}.
\end{equation*}
\end{punto}

\begin{proposition}
\label{mat-tens}\emph{(\cite[Theorem 3]{Kur2002})} Let $u$ and $v$ be
linearly recursive sequences over $R$ of orders $m,$ $n$ and with
characteristic polynomials $f(x),$ $g(x)$ respectively. Then

\begin{enumerate}
\item  $u\star _{g}v$ is a linearly recursive sequence over $R$ of order $%
m\cdot n$ and characteristic polynomial $\chi (S_{f}\otimes S_{g});$

\item  $u\star _{p}v$ is a linearly recursive sequence over $R$ of order $%
m\cdot n$ and characteristic polynomial $\chi (S_{f}\otimes
E_{n}+E_{m}\otimes S_{g}).$
\end{enumerate}
\end{proposition}

\begin{example}
Let $R$ be any ring and $\{x_{n}\}_{n=0}^{\infty },$ $\{y_{n}\}_{n=0}^{%
\infty }\in \mathcal{S}_{R}$ be solutions of the difference equations 
\begin{equation*}
\begin{tabular}{rrrrrr}
$x_{n+3}-x_{n+2}+x_{n-1}-x_{n}$ & $=$ & $0;$ & $x_{0}=0,$ & $x_{1}=1,$ & $%
x_{2}=2;$ \\ 
$y_{n+2}-y_{n+1}+y_{n}$ & $=$ & $0;$ & $y_{0}=1,$ & $y_{1}=0.$ & 
\end{tabular}
\end{equation*}
Then $\{x_{n}\}_{n=0}^{\infty }$ is a linearly recursive sequence over $R$
with characteristic polynomial $f(x)=x^{3}-x^{2}+x-1$ and $%
\{y_{n}\}_{n=0}^{\infty }$ is a linearly recursive sequence over $R$ with
characteristic polynomial $g(x)=x^{2}-x+1.$

Notice that 
\begin{eqnarray*}
S_{f}\otimes S_{g} &=&\left[ 
\begin{array}{ccc}
0 & 0 & 1 \\ 
1 & 0 & -1 \\ 
0 & 1 & 1
\end{array}
\right] \otimes \left[ 
\begin{array}{cc}
0 & -1 \\ 
1 & 1
\end{array}
\right] \\
&=&\left[ 
\begin{array}{cccccc}
0 & 0 & 0 & 0 & 0 & -1 \\ 
0 & 0 & 0 & 0 & 1 & 1 \\ 
0 & -1 & 0 & 0 & 0 & 1 \\ 
1 & 1 & 0 & 0 & -1 & -1 \\ 
0 & 0 & 0 & -1 & 0 & -1 \\ 
0 & 0 & 1 & 1 & 1 & 1
\end{array}
\right] .
\end{eqnarray*}
Hence $\{z_{n}\}_{n=0}^{\infty }:=\{x_{n}\}_{n=0}^{\infty }\star
_{g}\{y_{n}\}_{n=0}^{\infty }$ is by Proposition \ref{mat-tens} a linearly
recursive sequence over $R$ with characteristic polynomial 
\begin{equation*}
\chi (S_{f}\otimes S_{g})=x^{6}-x^{5}+x^{3}-x+1,
\end{equation*}
i.e. $\{z_{n}\}_{n=0}^{\infty }$ is a solution of the difference equation 
\begin{equation*}
z_{n+6}-z_{n+5}+z_{n+3}-z_{n+1}+z_{n}=0\text{ with initial vector }%
(0,0,-2,-1,0,1).
\end{equation*}
The following table gives the first $11$ terms of the sequences $%
\{z_{n}\}_{n=0}^{\infty }:$%
\begin{equation*}
\begin{tabular}{|c|c|c|c|c|c|c|c|c|c|c|c|}
\hline
$n$ & $0$ & $1$ & $2$ & $3$ & $4$ & $5$ & $6$ & $7$ & $8$ & $9$ & $10$ \\ 
\hline
$x_{n}$ & $0$ & $1$ & $2$ & $1$ & $0$ & $1$ & $2$ & $1$ & $0$ & $1$ & $2$ \\ 
\hline
$y_{n}$ & $1$ & $0$ & $-1$ & $-1$ & $0$ & $1$ & $1$ & $0$ & $-1$ & $-1$ & $0$
\\ \hline
$z_{n}$ & $0$ & $0$ & $-2$ & $-1$ & $0$ & $1$ & $2$ & $0$ & $0$ & $-1$ & $0$
\\ \hline
\end{tabular}
\end{equation*}
\end{example}

\begin{example}
Consider the sequences $\{x_{n}\}_{n=0}^{\infty }$ and $\{y_{n}\}_{n=0}^{%
\infty }$ of the previous example. Then 
\begin{equation*}
S_{f}\otimes E_{2}+E_{3}\otimes S_{g}=\left[ 
\begin{array}{cccccc}
0 & -1 & 0 & 0 & 1 & 0 \\ 
1 & 1 & 0 & 0 & 0 & 1 \\ 
1 & 0 & 0 & -1 & -1 & 0 \\ 
0 & 1 & 1 & 1 & 0 & -1 \\ 
0 & 0 & 1 & 0 & 1 & -1 \\ 
0 & 0 & 0 & 1 & 1 & 2
\end{array}
\right] .
\end{equation*}
By Proposition \ref{mat-tens} $\{z_{n}\}_{n=0}^{\infty
}=\{x_{n}\}_{n=0}^{\infty }\star _{p}\{y_{n}\}_{n=0}^{\infty
}:=\{\sum\limits_{j=0}^{n}\left( \QATOP{n}{j}\right) x_{j}\cdot
y_{n-j}\}_{n=0}^{\infty }$ is a linearly recursive sequence over $R$ with
characteristic polynomial 
\begin{equation*}
\chi (S_{f}\otimes E_{2}+E_{3}\otimes
S_{g})=x^{6}-5x^{5}+14x^{4}-25x^{3}+28x^{2}-15x+3.
\end{equation*}
Hence $\{z_{n}\}_{n=0}^{\infty }$ is a solution of the difference equation 
\begin{equation*}
z_{n+6}-5z_{n+5}+14z_{n+4}-25z_{n+3}+28z_{n+2}-15z_{n+1}+3z_{n}=0
\end{equation*}
with initial vector $(0,1,2,-2,-16,-29).$

The following table gives the first $9$ terms of the sequences $%
\{z_{n}\}_{n=0}^{\infty }:$ 
\begin{equation*}
\begin{tabular}{|c|c|c|c|c|c|c|c|c|c|}
\hline
$n$ & $0$ & $1$ & $2$ & $3$ & $4$ & $5$ & $6$ & $7$ & $8$ \\ \hline
$x_{n}$ & $0$ & $1$ & $2$ & $1$ & $0$ & $1$ & $2$ & $1$ & $0$ \\ \hline
$y_{n}$ & $1$ & $0$ & $-1$ & $-1$ & $0$ & $1$ & $1$ & $0$ & $-1$ \\ \hline
$z_{n}$ & $0$ & $1$ & $2$ & $-2$ & $-16$ & $-29$ & $-12$ & $29$ & $0$ \\ 
\hline
\end{tabular}
\end{equation*}
\end{example}

\begin{punto}
\textbf{Cofree comodules.}\label{ot C-adj} Let $C$ be an $R$-coalgebra. A
right $C$-comodule $(M,\varrho _{M})$ is called \emph{cofree}, if there
exists an $R$-module $K,$ such that $(M,\varrho _{M})\simeq (K\otimes
_{R}C,id_{K}\otimes \Delta _{C})$ as right $C$-comodules. Note that if $%
K\simeq R^{(\Lambda )},$ a free $R$-module, then $M\simeq R^{(\Lambda
)}\otimes _{R}C\simeq C^{(\Lambda )}$ as right $C$-comodules (this is one
reason of the terminology \emph{cofree}).
\end{punto}

As a direct consequence of Lemma \ref{cof} we get

\begin{corollary}
\label{LM*=}Let $M$ be an $R[\mathbf{x}]$-module. Then we have an
isomorphism of $R[\mathbf{x}]^{\circ }$-comodules 
\begin{equation*}
\mathcal{L}{\normalsize _{M^{\ast }}^{<k>}\simeq M[\mathbf{x}]^{\circ
}\simeq M^{\ast }\otimes _{R}R[\mathbf{x}]^{\circ }\simeq M^{\ast }\otimes
_{R}}\text{ }\mathcal{L}_{R}^{<k>}{\normalsize .}
\end{equation*}
In particular $M[\mathbf{x}]^{\circ }$ \emph{(}$\mathcal{L}_{M^{\ast
}}^{<k>} $\emph{)} is a cofree $R[\mathbf{x}]^{\circ }$-comodule \emph{(}$%
\mathcal{L}_{R}^{<k>}$-comodule\emph{)}.
\end{corollary}

\section{Linearly (bi)recursive bisequences}

\qquad In this section we consider the \emph{linearly }(\emph{bi})\emph{%
recursive }$k$\emph{-bisequences} and the \emph{reversible }$k$\emph{%
-sequences} over $R$-modules, where $R$ is an arbitrary commutative ground
ring. We generalize results of \cite{LT90} and \cite{KKMN95} concerning the
bialgebra structure of the linearly recursive sequences over a base field to
the case of arbitrary \emph{artinian} ground rings.

\begin{punto}
Let $M$ be an $R$-module, $\mathbf{l=}(l_{1},...,l_{k})\in \mathbb{N}%
_{0}^{k} $ and consider the \emph{system of linear bidifference equations }%
(ab. SLBE) 
\begin{equation}
\begin{tabular}{lllll}
$x_{\mathbf{z}+(l_{1},0,...,0)}$ & $+$ & $\sum%
\limits_{i=1}^{l_{1}}p_{(1,l_{1}-i)}(\mathbf{z})x_{\mathbf{z}%
+(l_{1}-i,0,...,0)}$ & $=$ & $g_{1}(\mathbf{z}),$ \\ 
$x_{\mathbf{z}+(0,l_{2},0,...,0)}$ & $+$ & $\sum%
\limits_{i=1}^{l_{2}}p_{(2,l_{2}-i)}(\mathbf{z})x_{\mathbf{z}%
+(0,l_{2}-i,0,...,0)}$ & $=$ & $g_{2}(\mathbf{z}),$ \\ 
$...$ & $...$ & $...$ & $...$ & $...$ \\ 
$...$ & $...$ & $...$ & $...$ & $...$ \\ 
$x_{\mathbf{z}+(0,...,0,l_{k})}$ & $+$ & $\sum%
\limits_{i=1}^{l_{k}}p_{(k,l_{k}-i)}(\mathbf{z})x_{\mathbf{z}%
+(0,...,0,l_{k}-i)}$ & $=$ & $g_{k}(\mathbf{z}),$%
\end{tabular}
\label{SLBE}
\end{equation}
where the $p_{jl}$'s are $R$-valued functions and the $g_{j}$'s are $M$%
-valued functions defined for all $\mathbf{z}\in \mathbb{Z}^{<k>}.$ If the $%
g_{j}$'s are identically zero, then (\ref{SLBE}) is said to be a \emph{%
homogenous }SLBE.\emph{\ }If the $p_{jl}$'s are constants, then (\ref{SLBE})
is said to be a SLBE \emph{with constant coefficients}.
\end{punto}

\begin{punto}
\textbf{Bisequences.}\label{bis} For an $R$-module $M$ and $k\geq 0$ let 
\begin{equation*}
\mathbb{\widetilde{\mathcal{S}}}_{M}^{<k>}:=\{\widetilde{\nu }:\mathbb{Z}%
^{k}\rightarrow M\}\simeq M^{\mathbb{Z}^{k}}
\end{equation*}
be the $R$-module of $k$-\emph{bisequences} over $M.$ If $M$ (resp. $k$) is
not mentioned, then we mean $M=R$ (resp. $k=1$). For $\widetilde{w}\in 
\mathbb{\widetilde{\mathcal{S}}}_{M}^{<k>}$ and $f(\mathbf{x})=\sum\limits_{%
\mathbf{i}}a_{\mathbf{i}}\mathbf{x}^{\mathbf{i}}\in R[\mathbf{x},\mathbf{x}%
^{-1}]$ define 
\begin{equation*}
f(\mathbf{x})\rightharpoonup \widetilde{w}=\widetilde{\nu }\in \mathbb{%
\widetilde{\mathcal{S}}}_{M}^{<k>},\text{ where }\widetilde{\nu }(\mathbf{z}%
):=\sum_{\mathbf{i}}a_{\mathbf{i}}\widetilde{w}(\mathbf{z}+\mathbf{i})\text{
for all }\mathbf{z}\in \mathbb{Z}^{k}.
\end{equation*}
With this action $\mathbb{\widetilde{\mathcal{S}}}_{M}^{<k>}$ becomes an $R[%
\mathbf{x},\mathbf{x}^{-1}]$-module. For subsets $I\subset R[\mathbf{x},%
\mathbf{x}^{-1}]$ and $Y\subset \mathbb{\widetilde{\mathcal{S}}}_{M}^{<k>}$
consider 
\begin{equation*}
\begin{tabular}{lll}
$\mathrm{An}_{\mathbb{\widetilde{\mathcal{S}}}_{M}^{<k>}}(I)$ & $=$ & $\{%
\widetilde{w}\in \mathbb{\widetilde{\mathcal{S}}}_{M}^{<k>}\,|$ $%
g\rightharpoonup \widetilde{w}=0$ for every $g\in I\},$ \\ 
$\mathrm{An}_{R[\mathbf{x},\mathbf{x}^{-1}]}(Y)$ & $=$ & $\{h\in R[\mathbf{x}%
,\mathbf{x}^{-1}]\,|\,h\rightharpoonup \widetilde{\nu }=0$ for every $%
\widetilde{\nu }\in Y\}.$%
\end{tabular}
\end{equation*}
Obviously $\mathrm{An}_{\mathcal{S}_{M}}^{<k>}(I)\subset \mathcal{S}%
_{M}^{<k>}$ is an $R[\mathbf{x},\mathbf{x}^{-1}]$-submodule and $\mathrm{An}%
_{R[\mathbf{x},\mathbf{x}^{-1}]}(Y)\vartriangleleft $ $R[\mathbf{x},\mathbf{x%
}^{-1}]\,$ is an ideal.
\end{punto}

\begin{definition}
Let $M$ be an $R$-module. We call $\widetilde{w}\in \mathbb{\widetilde{%
\mathcal{S}}}_{M}^{<k>}$ a \emph{linearly recursive }$k$\emph{-bisequence }
(resp. a \emph{linearly birecursive }$k$\emph{-bisequence}), if $\mathrm{An}%
_{R[\mathbf{x}]}(\widetilde{w})$ is a monic ideal (resp. a reversible
ideal). Note that a $k$-bisequence $\widetilde{u}\in \mathbb{\widetilde{%
\mathcal{S}}}_{M}^{<k>}$ is linearly recursive,\emph{\ }iff it's a solution
of a homogenous SLBE with constants coefficients of the form (\ref{SLBE}).
The subsets $\widetilde{\mathcal{L}}_{M}^{<k>}\subseteq \mathbb{\widetilde{%
\mathcal{S}}}_{M}^{<k>}$ of linearly recursive $k$-bisequences and $%
\widetilde{\mathcal{B}}_{M}^{<k>}\subseteq \mathbb{\widetilde{\mathcal{S}}}%
_{M}^{<k>}$ of linearly birecursive $k$-bisequences over $M$ are obviously $%
R[\mathbf{x},\mathbf{x}^{-1}]$-submodules.
\end{definition}

\section*{Reversible sequences over modules}

\begin{punto}
\label{rever}Let $M$ be an $R$-module. A $k$-bisequence $\widetilde{u}$ is
said to be a \emph{reverse} of $u\in \mathcal{S}_{M}^{<k>}$, if $\widetilde{u%
}|_{\mathbb{N}_{0}^{k}}=u$\ and $\mathrm{An}_{R[\mathbf{x}]}(\widetilde{u})=%
\mathrm{An}_{R[\mathbf{x}]}(u).$ A linearly recursive $k$-sequence $u$ will
be called \emph{reversible}, if $u$ has a reverse $\widetilde{u}\in 
\widetilde{\mathcal{L}}_{M}^{<k>}.$ With $\mathcal{R}_{M}^{<k>}\subset 
\mathcal{L}_{M}^{<k>}$ we denote the $R[\mathbf{x}]$-submodule of reversible 
$k$-sequences over $M.$
\end{punto}

\begin{lemma}
\label{qx}\emph{(}Compare\emph{\ \cite[Proposition 14.11]{KKMN95})} Let $R$
be \emph{artinian}.

\begin{enumerate}
\item  Every monic ideal $I\vartriangleleft R[\mathbf{x}]$ contains a subset
of monic polynomials 
\begin{equation}
\{x_{j}^{d_{j}}q_{j}(x_{j})|\text{ }q_{j}(x_{j})\text{ is reversible for }%
j=1,...,k\}.  \label{revers}
\end{equation}

\item  Let $M$ be an $R$-module. Then every linearly recursive $k$%
-bisequence over $M$ is linearly birecursive \emph{(}i.e. $\widetilde{%
\mathcal{B}}_{M}^{<k>}=\widetilde{\mathcal{L}}_{M}^{<k>}$\emph{)}.
\end{enumerate}
\end{lemma}

\begin{Beweis}
\begin{enumerate}
\item  By \cite[8.7]{AM69} every commutative artinian ring is (up to
isomorphism)\ a direct sum of local artinian rings. W.l.o.g. let $R$ be a
local artinian ring. The Jacobson radical of $R$%
\begin{equation*}
J(R)=\{r\in R|\text{ }r\text{ is not invertible in }R\}
\end{equation*}
is nilpotent, hence there exists a positive integer $n,$ such that $%
J(R)^{n}=0.$ Let $I$ be a monic ideal with a subset of monic polynomials $%
\{g_{1}(x_{1}),...,g_{k}(x_{k})\}\subset I.$ If $g_{j}(x_{j})\equiv
f_{j}(x_{j})$ $(\mathrm{mod}$ $J(R)[x_{j}])$ for $j=1,...,k,$ then $%
g_{j}(x_{j})|f_{j}(x_{j})^{n},$ where $n$ is the index of nilpotency of the
ideal $J(R).$ Hence $f_{j}(x_{j})^{n}\in I.$ If we write $%
f_{j}(x_{j})^{n}=x_{j}^{d_{j}}q_{j}(x_{j})$ with $(x_{j},q_{j}(x_{j}))=1,$
then $q_{j}(0)\in U(R),$ i.e. $q_{j}(x_{j})$ is a reversible polynomial for $%
j=1,...,k.$

\item  Let $\widetilde{u}$ be a linearly recursive $k$-bisequence over $M.$
If $R$ is artinian, then $\mathrm{An}_{R[\mathbf{x}]}(\widetilde{u})$
contains by (1) a subset of monic polynomials $\{x_{j}^{d_{j}}q_{j}(x_{j})|$ 
$q_{j}(x_{j})$ is reversible for $j=1,...,k\}.$ Then for every $\mathbf{z}%
\in \mathbb{Z}^{k}$ we have $(q_{j}(x_{j})\rightharpoonup \widetilde{u}%
)(z_{1},...,z_{j},...,z_{k})=(x_{j}^{d_{j}}q_{j}(x_{j})\rightharpoonup 
\widetilde{u})(z_{1},...,z_{j}-d_{j},...,z_{k})=0.$ Hence $\{q_{j}(x_{j})|$ $%
i=1,...,k\}\subset \mathrm{An}_{R[\mathbf{x}]}(\widetilde{u}),$ i.e. $%
\mathrm{An}_{R[\mathbf{x}]}(\widetilde{u})$ is a reversible ideal.$%
\blacksquare $
\end{enumerate}
\end{Beweis}

\begin{punto}
\textbf{Backsolving}.\label{BS} Let $M$ be an $R$-module. Let $u$ be a
linearly recursive sequence over $M$ and assume that $\mathrm{An}_{R[x]}(u)$
contains some monic polynomial of the form $%
x^{d}q(x)=x^{d}(a_{0}+a_{1}x+...+a_{l-1}x^{l-1}+x^{l}),$ $a_{0}\in U(R).$
Then 
\begin{equation*}
a_{0}u(j+d)+a_{1}u(j+d+1)+...+a_{l-1}u(j+d+l-1)+u(j+d+l)=0\text{ for all }%
j\geq 0
\end{equation*}
and we get by \emph{Backsolving} a \emph{unique }linearly birecursive
bisequence $\widetilde{u}\in \mathrm{An}_{\widetilde{\mathcal{S}}_{M}}(q(x))$
with $\widetilde{u}(n)=u(n)$ for all $n\geq d$. The bisequence $\widetilde{u}%
\equiv 0$ in case $l=0$ and is given for $l\neq 0$ by 
\begin{equation*}
\widetilde{u}(z):= 
\begin{cases}
u(z), & z\geq d \\ 
-a_{0}^{-1}(a_{1}\widetilde{u}(z+1)+...+a_{l-1}\widetilde{u}(z+l-1)+%
\widetilde{u}(z+l)), & z<d.
\end{cases}
\end{equation*}
If there are two bisequences $\widetilde{v},$ $\widetilde{w}\in \mathrm{An}_{%
\widetilde{\mathcal{S}}_{M}}(q(x))$ with $\widetilde{v}(n)=u(n)=\widetilde{w}%
(n)$ for all $n\geq d,$ then one can easily show by backsolving using $q(x)$
that $\widetilde{v}=\widetilde{w}.$ Moreover we claim that $\mathrm{An}%
_{R[x]}(\widetilde{u})=\mathrm{An}_{R[x]}(u).$ It's obvious that $\mathrm{An}%
_{R[x]}(\widetilde{u})\subseteq \mathrm{An}_{R[x]}(u).$ On the other hand
assume $g(x)=\sum\limits_{j=0}^{m}b_{j}x^{j}\in \mathrm{An}_{R[x]}(u).$ We
prove by induction that $(g\rightharpoonup \widetilde{u})(z)=0$ for all $%
z\in \mathbb{Z}.$ First of all, note that for all $z\geq d$ we have $%
(g\rightharpoonup \widetilde{u})(z)=(g\rightharpoonup u)(z)=0.$ Now let $%
z_{0}<d$ and assume that $(g\rightharpoonup \widetilde{u})(z)=0$ for $z\in
\{z_{0},z_{0}+1,...,z_{0}+l-1\}\subseteq \mathbb{Z}.$ Then we have for $%
z=z_{0}-1:$%
\begin{equation*}
\begin{tabular}{lll}
$(g\rightharpoonup \widetilde{u})(z_{0}-1)$ & $=$ & $\sum%
\limits_{j=0}^{m}b_{j}\widetilde{u}(j+z_{0}-1)$ \\ 
& $=$ & $\sum\limits_{j=0}^{m}b_{j}(\sum\limits_{i=1}^{l}-a_{0}^{-1}a_{i}%
\widetilde{u}(j+z_{0}-1+i))$ \\ 
& $=$ & $-\sum\limits_{i=1}^{l}a_{0}^{-1}a_{i}\sum\limits_{j=0}^{m}b_{j}%
\widetilde{u}(j+z_{0}-1+i)$ \\ 
& $=$ & $-\sum\limits_{i=1}^{l}a_{0}^{-1}a_{i}(g\rightharpoonup \widetilde{u}%
)(z_{0}-1+i)$ \\ 
& $=$ & $0.$%
\end{tabular}
\end{equation*}

If $u$ is a linearly recursive $k$-sequence over $M$ with $k>1$ and $\mathrm{%
An}_{R[\mathbf{x}]}(u)$ contains a set of monic polynomials $%
\{x_{j}^{d_{j}}q_{j}(x_{j})\mid q_{j}$ is reversible for $j=1,...,k\},$ then
we get by \emph{backsolving} through $q_{j}(x_{j})$ along the $j$-th row for 
$j=1,...,k$ a unique linearly birecursive $k$-bisequence $\widetilde{u}\in 
\mathrm{An}_{\mathbb{\widetilde{\mathcal{S}}}%
_{M}^{<k>}}(q_{1}(x_{1}),...,q_{k}(x_{k}))$ with $\widetilde{u}(\mathbf{n}%
)=u(\mathbf{n})$ for all $\mathbf{n}\geq \mathbf{d}$ and it follows moreover
that $\mathrm{An}_{R[\mathbf{x}]}(\widetilde{u})=\mathrm{An}_{R[\mathbf{x}%
]}(u).$
\end{punto}

\begin{lemma}
\label{s-rev}Let $M$ be an $R$-module.

\begin{enumerate}
\item  Every birecursive $k$-sequence over $M$ is reversible with unique
reverse \emph{(}which we denote by $\mathrm{Rev}(u)$\emph{)}. Moreover $%
\mathcal{B}_{M}^{<k>}$ becomes a structure of an $R[\mathbf{x},\mathbf{x}%
^{-1}]$-module through $f\rightharpoonup u:=(f\rightharpoonup \mathrm{Rev}%
(u))_{\mid _{\mathbb{N}_{0}^{k}}}.$

\item  If $R$ is artinian, then every reversible $k$-sequence over $M$ is
birecursive as well \emph{(}i.e. $\mathcal{B}_{M}^{<k>}=\mathcal{R}%
_{M}^{<k>} $\emph{)}.
\end{enumerate}
\end{lemma}

\begin{Beweis}
Let $M$ be an $R$-module.

\begin{enumerate}
\item  If $u\in \mathcal{B}_{M}^{<k>},$ then $\mathrm{An}_{R[\mathbf{x}]}(u)$
contains a set of reversible polynomials $\{q_{j}(x_{j})\mid j=1,...,k\}$
and we get by backsolving (see \ref{BS}) a \emph{unique }linearly
birecursive $k$-bisequence $\widetilde{u}\in \mathrm{An}_{\mathbb{\widetilde{%
\mathcal{S}}}_{M}^{<k>}}(q_{1}(x_{1}),...,q_{k}(x_{k}))$ with $\widetilde{u}(%
\mathbf{n})=u(\mathbf{n})$ for all $\mathbf{n}\in \mathbb{N}_{0}^{k}.$ For
the bisequence $\widetilde{u}$ we have as shown above $\mathrm{An}_{R[%
\mathbf{x}]}(\widetilde{u})=\mathrm{An}_{R[\mathbf{x}]}(u),$ i.e. $%
\widetilde{u}$ is a reverse of $u.$ The last statement is obvious.

\item  By (1) $\mathcal{B}_{M}^{<k>}\subseteq \mathcal{R}_{M}^{<k>}.$ If $R$
is artinian and $u\in \mathcal{R}_{M}^{<k>}$ with reverse $\widetilde{u},$
then $\mathrm{An}_{R[\mathbf{x}]}(u)=\mathrm{An}_{R[\mathbf{x}]}(\widetilde{u%
})$ is by Lemma \ref{qx} (2) reversible, i.e. $u\in \mathcal{B}%
_{M}^{<k>}.\blacksquare $
\end{enumerate}
\end{Beweis}

\begin{example}
The Fibonacci sequence $\digamma =(0,1,1,2,3,5,...)$ has elementary
characteristic polynomial $f(x)=x^{2}-x-1.$ Since $f(0)=-1$ is invertible in 
$\mathbb{Z},$ we conclude that $\digamma $ is reversible with reverse 
\begin{equation*}
\mathrm{Rev}(\digamma )(z)= 
\begin{cases}
\digamma (z) & z\geq 0 \\ 
\mathrm{Rev}(\digamma )(z+2)-\mathrm{Rev}(\digamma )(z+1) & z<0.
\end{cases}
\end{equation*}
The following tables lists some of the terms of the bisequence $\mathrm{Rev}%
(\digamma )\in \mathrm{An}_{\mathcal{S}_{\mathbb{Z}}}(x^{2}-x-1):$%
\begin{equation*}
\begin{tabular}{|l|l|l|l|l|l|l|l|l|l|l|}
\hline
$z$ & $...$ & $-4$ & $-3$ & $-2$ & $-1$ & $0$ & $1$ & $2$ & $3$ & $4$ \\ 
\hline
$\mathrm{Rev}(\digamma )(z)$ & $...$ & $-3$ & $2$ & $-1$ & $1$ & $0$ & $1$ & 
$1$ & $2$ & $3$ \\ \hline
\end{tabular}
\end{equation*}
\end{example}

\begin{lemma}
\label{bi=rev}We have an isomorphism of $R[\mathbf{x},\mathbf{x}^{-1}]$%
-modules 
\begin{equation}
\widetilde{\mathcal{B}}_{M}^{<k>}\simeq \mathcal{B}_{M}^{<k>}.
\label{B=telB}
\end{equation}
\end{lemma}

\begin{Beweis}
By Lemma \ref{s-rev} we have the well defined $R[\mathbf{x},\mathbf{x}^{-1}]$%
-linear map 
\begin{equation*}
\mathrm{Rev}(-):\mathcal{B}_{M}^{<k>}\rightarrow \mathbb{\widetilde{\mathcal{%
B}}}_{M}^{<k>},\text{ }u\mapsto \mathrm{Rev}(u).
\end{equation*}
It's easy to see that $\mathrm{Rev}(-)$ is bijective with inverse $%
\widetilde{u}\mapsto \widetilde{u}_{\mid _{\mathbb{N}_{0}^{k}}}.\blacksquare 
$
\end{Beweis}

\begin{punto}
\label{per-deg}Let $M$ be an $R$-module. We call a $k$-sequence $u\in 
\mathcal{S}_{M}^{<k>}$ \emph{periodic} (resp. \emph{degenerating}),\emph{\ }%
if $\mathbf{x}^{\mathbf{d}}(\mathbf{x}^{\mathbf{t}}\rightharpoonup u)=0$ for
some $\mathbf{d}\in \mathbb{N}_{0}^{k}$ and $\mathbf{t}\in \mathbb{N}^{k}$
(resp. $\mathbf{x}^{\mathbf{d}}\rightharpoonup u=0$ for some $\mathbf{d}\in 
\mathbb{N}_{0}^{k}$). It's clear that the subsets $\mathcal{P}%
_{M}^{<k>}\subseteq \mathcal{L}_{M}^{<k>}$ of periodic $k$-sequences and $%
\mathcal{D}_{M}^{<k>}\subseteq \mathcal{L}_{M}^{<k>}$ of degenerating $k$%
-sequences are $R[\mathbf{x}]$-submodules.
\end{punto}

\begin{remark}
(\cite[Proposition 5.2]{KKMN95}) If $M$ is a \emph{finite} $R$-module, then
every linearly recursive sequence over $M$ is periodic (i.e. $\mathcal{P}%
_{M}^{<1>}=\mathcal{L}_{M}^{<1>}$).
\end{remark}

\begin{proposition}
\label{P=D+R}\emph{(\cite[Proposition 5.27]{KKMN95}) }Let $R$ be an \emph{%
arbitrary} commutative ring, $M$ an $R$-module and denote with $\mathcal{R}%
\mathcal{P}_{M}^{<k>}$ the set of reversible periodic $k$-sequences over $M.$
Then we have an isomorphism of $R[\mathbf{x}]$-modules 
\begin{equation}
\mathcal{P}_{M}^{<k>}\simeq \mathcal{D}_{M}^{<k>}\oplus \mathcal{RP}%
_{M}^{<k>}.  \label{per=}
\end{equation}
\end{proposition}

The following result generalizes Proposition \ref{P=D+R} and describes the $%
R[\mathbf{x}]$-module structure of arbitrary linearly recursive $k$%
-sequences of $R$-modules, where $R$ is an \emph{artinian} commutative
ground ring:

\begin{proposition}
\label{split}Let $M$ be an $R$-module. If $R$ is artinian, then we have
isomorphisms of $R[\mathbf{x}]$-modules 
\begin{equation}
{\normalsize \mathcal{L}}_{M}^{<k>}\simeq \mathcal{D}_{M}^{<k>}\oplus 
\mathbb{\widetilde{{\normalsize \mathcal{L}}}}_{M}^{<k>}=\mathcal{D}%
_{M}^{<k>}\oplus \widetilde{\mathcal{B}}_{M}^{<k>}\simeq \mathcal{D}%
_{M}^{<k>}\oplus \mathcal{B}_{M}^{<k>}=\mathcal{D}_{M}^{<k>}\oplus \mathcal{R%
}_{M}^{<k>}.  \label{LM=}
\end{equation}
\end{proposition}

\begin{Beweis}
Let $R$ be artinian and $M$ an $R$-module. If $u$ is a linearly recursive
sequence over $M,$ then $\mathrm{An}_{R[\mathbf{x}]}(u)$ contains by Lemma 
\ref{qx} (1) a set of monic polynomials $\{x^{d_{j}}q_{j}(x_{j})\mid
q_{j}(x_{j})$ is reversible for $j=1,...,k\}.$ By \emph{backsolving }(see 
\ref{BS}) we have a well defined morphism of $R[\mathbf{x}]$-modules 
\begin{equation}
\gamma :{\normalsize \mathcal{L}}_{M}^{<k>}\rightarrow \mathbb{\widetilde{%
{\normalsize \mathcal{L}}}}_{M}^{<k>},\text{ }u\mapsto \widetilde{u},
\label{gama}
\end{equation}
where $\widetilde{u}$ is the \emph{unique} linearly birecursive bisequence $%
\widetilde{u}\in \mathrm{An}_{\widetilde{\mathcal{S}}_{M}}(q_{1},...,q_{k})$
with $\widetilde{u}(\mathbf{n})=u(\mathbf{n})$ for all $\mathbf{n}\geq 
\mathbf{d}.$ It's clear that $\mathrm{Ke}(\gamma )=\mathcal{D}_{M}^{<k>}.$
On the other hand, there is a morphism of $R[\mathbf{x}]$-modules 
\begin{equation}
\beta ={\normalsize \widetilde{\mathcal{L}}}_{M}^{<k>}\rightarrow \mathbb{%
{\normalsize \mathcal{L}}}_{M}^{<k>},\text{ }\widetilde{w}\mapsto \widetilde{%
w}_{|_{\mathbb{N}_{0}^{k}}}.  \label{beta}
\end{equation}
It's obvious that $\gamma \circ \beta =id_{\widetilde{\mathcal{L}}%
_{M}^{<k>}},$ hence the following exact sequence of $R[\mathbf{x}]$-modules
splits 
\begin{equation*}
0\rightarrow \mathcal{D}_{M}^{<k>}\rightarrow {\normalsize \mathcal{L}}%
_{M}^{<k>}\overset{\gamma }{\rightarrow }\mathbb{\widetilde{{\normalsize 
\mathcal{L}}}}_{M}^{<k>}\rightarrow 0,
\end{equation*}
i.e. $\mathcal{L}_{M}^{<k>}\simeq \mathcal{D}_{M}^{<k>}\oplus \mathbb{%
\widetilde{\mathcal{L}}}_{M}^{<k>}.$ Since $R$ is artinian, we have by
Lemmata \ref{qx} (2) and \ref{s-rev} (2) $\mathbb{\widetilde{\mathcal{L}}}%
_{M}^{<k>}=\mathbb{\widetilde{\mathcal{B}}}_{M}^{<k>}$ and $\mathbb{\mathcal{%
B}}_{M}^{<k>}=\mathcal{R}_{M}^{<k>}.$ We are done now by the isomorphism of $%
R[\mathbf{x}]$-modules $\mathcal{B}_{M}^{<k>}\simeq \mathbb{\widetilde{%
\mathcal{B}}}_{M}^{<k>}$ (Lemma \ref{bi=rev}).$\blacksquare $
\end{Beweis}

\begin{punto}
\textbf{The Hopf }$R$-\textbf{algebra }$R[\mathbf{x},\mathbf{x}^{-1}].$
Consider the commutative group $G$ generated by $\{x_{j}\mid j=1,...,k\}.$
Then the ring of Laurent polynomials $R[\mathbf{x},\mathbf{x}^{-1}]=RG$ has
the structure of a \emph{commutative cocommutative} Hopf $R$-algebra $(R[%
\mathbf{x},\mathbf{x}^{-1}],\mu ,\eta ,\Delta ,\varepsilon ,S),$ where $\mu $
resp. $\eta $ are the usual multiplication resp. the usual unity and 
\begin{equation}
\begin{tabular}{lllllll}
$\Delta :$ & $R[\mathbf{x},\mathbf{x}^{-1}]$ & $\rightarrow $ & $R[\mathbf{x}%
,\mathbf{x}^{-1}]\otimes _{R}R[\mathbf{x},\mathbf{x}^{-1}],$ & $%
x_{j}^{z}\mapsto $ & $x_{j}^{z}\otimes x_{j}^{z},$ & $\forall $ $z\in 
\mathbb{Z},$ $j=1,...,k,$ \\ 
$\varepsilon :$ & $R[\mathbf{x},\mathbf{x}^{-1}]$ & $\rightarrow $ & $R,$ & $%
x_{j}^{z}\mapsto $ & $1_{R},$ & $\forall $ $z\in \mathbb{Z},$ $j=1,...,k,$
\\ 
$S:$ & $R[\mathbf{x},\mathbf{x}^{-1}]$ & $\rightarrow $ & $R[\mathbf{x},%
\mathbf{x}^{-1}],$ & $x_{j}^{z}\mapsto $ & $x_{j}^{-z},$ & $\forall $ $z\in 
\mathbb{Z},$ $j=1,...,k.$%
\end{tabular}
\label{Lk-coal}
\end{equation}
\end{punto}

\begin{proposition}
\label{rx-1-admiss}Let $R$ be an \emph{arbitrary} commutative ring. Then $R[%
\mathbf{x},\mathbf{x}^{-1}]$ is an admissible Hopf $R$-algebra and $R[%
\mathbf{x},\mathbf{x}^{-1}]^{\circ }$ is a Hopf $R$-algebra.
\end{proposition}

\begin{Beweis}
Notice that $R[\mathbf{x},\mathbf{x}^{-1}]$ is a cofinitary Hopf $R$-algebra
by Lemma \ref{cof} (2). Consider the proof of Proposition \ref{rx-admiss}
and replace $R[\mathbf{x}]$ with $R[\mathbf{x},\mathbf{x}^{-1}].$ Then the
map 
\begin{equation*}
T_{j}:B\rightarrow B,\text{ }b\mapsto \overline{\Delta }(x_{j})b
\end{equation*}
is invertible with inverse 
\begin{equation*}
\overline{T}_{j}:B\rightarrow B,\text{ }b\mapsto \overline{\Delta }%
(x_{j}^{-1})b.
\end{equation*}
Then the matrix $M_{j}$ of $T_{j}$ is invertible and $\chi _{j}(0)\in U(R)$
for $j=1,...,k.$ Consequently $\mathcal{K}_{R[\mathbf{x},\mathbf{x}^{-1}]}$
satisfies axiom (\ref{A1}). Since $R[\mathbf{x},\mathbf{x}^{-1}]/\mathrm{Ke}%
(\varepsilon )\simeq R,$ $\mathcal{K}_{R[\mathbf{x};p]}$ satisfies axiom (%
\ref{A2}). Consider the \emph{bijective }antipode $S$ of $R[\mathbf{x};p].$
For every ideal $I\vartriangleleft R[\mathbf{x},\mathbf{x}^{-1}],$ $%
S^{-1}(I)\vartriangleleft R[\mathbf{x},\mathbf{x}^{-1}]$ is an ideal and we
have an isomorphism of $R$-modules $R[\mathbf{x},\mathbf{x}%
^{-1}]/S^{-1}(I)\simeq R[\mathbf{x},\mathbf{x}^{-1}]/I.$ Hence $\mathcal{K}%
_{R[\mathbf{x},\mathbf{x}^{-1}]}$ satisfies axiom (\ref{A3}). Consequently $%
R[\mathbf{x},\mathbf{x}^{-1}]$ is an admissible Hopf $R$-algebra. The last
statement follows now by Proposition \ref{cofinitary}.$\blacksquare $
\end{Beweis}

For every $R$-module $M$ we have an isomorphism of $R[\mathbf{x},\mathbf{x}%
^{-1}]$-modules 
\begin{equation}
\Psi _{M}:M[\mathbf{x},\mathbf{x}^{-1}]^{\ast }\rightarrow \mathbb{%
\widetilde{\mathcal{S}}}_{M^{\ast }}^{<k>}\mathbb{,}\text{ }\widetilde{%
\varphi }\mapsto \lbrack \mathbf{z}\mapsto \lbrack m\mapsto \widetilde{%
\varphi }(m\mathbf{x}^{\mathbf{z}})]]  \label{Psi_M}
\end{equation}
with inverse $\widetilde{u}\mapsto \lbrack m\mathbf{x}^{\mathbf{z}}\mapsto 
\widetilde{u}(\mathbf{z})(m)].$

As in the proof of Proposition \ref{mxn-d-l} we get

\begin{proposition}
\label{corr}Let $R$ be an arbitrary ring. Then \emph{(\ref{Psi_M}) }induces
an isomorphism of $R[\mathbf{x},\mathbf{x}^{-1}]$-modules 
\begin{equation}
M[\mathbf{x},\mathbf{x}^{-1}]^{\circ }\simeq \widetilde{\mathcal{B}}%
_{M^{\ast }}^{<k>}.  \label{Lau=srev}
\end{equation}
\end{proposition}

\begin{Beweis}
Consider the isomorphism of $R[\mathbf{x},\mathbf{x}^{-1}]$-modules $M[%
\mathbf{x},\mathbf{x}^{-1}]^{\ast }\overset{\Psi _{M}}{\simeq }\mathbb{%
\widetilde{\mathcal{S}}}_{M^{\ast }}^{<k>}$ (\ref{Psi_M}). Let $\varkappa
\in M[\mathbf{x},\mathbf{x}^{-1}]^{\circ }.$ Then $I\rightharpoonup
\varkappa =0$ for some $R$-cofinite $R[\mathbf{x},\mathbf{x}^{-1}]$-ideal $%
I\vartriangleleft R[\mathbf{x},\mathbf{x}^{-1}]$ and so $I\rightharpoonup
\Psi (\varkappa )=\Psi (I\rightharpoonup \varkappa )=0.$ By Lemma \ref{cof}
(2) $I$ is a reversible ideal and so $\mathrm{An}_{R[\mathbf{x}]}(u)\supset
I\cap R[\mathbf{x}]$ is a reversible ideal, i.e. $\Psi (\varkappa )$ is
linearly birecursive.

On the other hand, let $\widetilde{u}\in \widetilde{\mathcal{B}}_{M^{\ast
}}^{<k>}.$ Then $\mathrm{An}_{R[\mathbf{x}]}(\widetilde{u})$ is by
definition a reversible ideal, i.e. it contains a subset of reversible
polynomials $\{q_{j}(x_{j}),$ $j=1,...,k\}.$ Note that for arbitrary $g\in R[%
\mathbf{x},\mathbf{x}^{-1}]$ we have $gq_{j}\rightharpoonup \Psi ^{-1}(%
\widetilde{u})=\Psi ^{-1}(gq_{j}\rightharpoonup \widetilde{u})=\Psi
^{-1}(g\rightharpoonup (q_{j}\rightharpoonup \widetilde{u}))=0$ for $%
j=1,...,k.$ By Lemma \ref{cof} (2) the reversible ideal $%
(q_{1}(x_{1}),...,q_{k}(x_{k}))\vartriangleleft R[\mathbf{x},\mathbf{x}%
^{-1}] $ is $R$-cofinite, i.e. $\Psi ^{-1}(\widetilde{u})\in M[\mathbf{x},%
\mathbf{x}^{-1}]^{\circ }.\blacksquare $
\end{Beweis}

\begin{punto}
\textbf{The Hopf }$R$-\textbf{algebra structures on }$\widetilde{\mathcal{B}}%
^{<k>}$ and $\mathcal{B}^{<k>}.$\label{la-co} Let $R$ be an arbitrary ring
and consider the Hopf $R$-algebra $R[\mathbf{x},\mathbf{x}^{-1}].$ Then $%
\widetilde{\mathcal{S}}^{<k>}\simeq R^{\mathbb{Z}^{k}}\simeq R[\mathbf{x},%
\mathbf{x}^{-1}]^{\ast }$ is an $R$-algebra with the \emph{Hadamard product} 
\begin{equation}
\star :{\normalsize \widetilde{\mathcal{S}}^{<k>}\otimes _{R}\widetilde{%
\mathcal{S}}^{<k>}\rightarrow \widetilde{\mathcal{S}}^{<k>},}\text{ }%
\widetilde{u}\otimes \widetilde{v}\mapsto \lbrack \mathbf{z}\mapsto 
\widetilde{u}(\mathbf{z})\widetilde{v}(\mathbf{z})]  \label{Had-2}
\end{equation}
and the unity 
\begin{equation}
\eta :R\rightarrow \widetilde{\mathcal{S}}{\normalsize ^{<k>}},\text{ }%
1_{R}\mapsto \lbrack \mathbf{z}\mapsto 1_{R}]\text{ for every }\mathbf{z}\in 
\mathbb{Z}^{k}.  \label{et-bi}
\end{equation}
By Proposition \ref{rx-1-admiss} $R[\mathbf{x},\mathbf{x}^{-1}]^{\circ }$ is
a Hopf $R$-algebra. So $\mathcal{B}^{<k>}\simeq R[\mathbf{x},\mathbf{x}%
^{-1}]^{\circ }$ inherits the structure of a Hopf $R$-algebra $(\mathcal{B}%
^{<k>},\star _{g},\eta _{g},\Delta _{\mathcal{B}^{<k>}},\varepsilon _{%
\mathcal{B}^{<k>}},S_{\mathcal{B}^{<k>}}),$ where $\star _{g}$ is the
Hadamard product (\ref{Hadamard}), $\eta _{g}$ is the unity (\ref{eta-g})
and 
\begin{equation*}
\begin{tabular}{llllllll}
$\Delta _{\mathcal{B}^{<k>}}$ & $:$ & $\mathcal{B}^{<k>}$ & $\rightarrow $ & 
$\mathcal{B}^{<k>}\otimes _{R}\mathcal{B}^{<k>},$ & $u$ & $\mapsto $ & $%
\sum\limits_{\mathbf{t}\leq \mathbf{l}-\mathbf{1}}(\mathbf{x}^{\mathbf{t}%
}\rightharpoonup u)\otimes e_{\mathbf{t}}^{\mathbf{F}},$ \\ 
$\varepsilon _{\mathcal{B}^{<k>}}$ & $:$ & $\mathcal{B}^{<k>}$ & $%
\rightarrow $ & $R,$ & $u$ & $\mapsto $ & $u(\mathbf{0}),$ \\ 
$S_{\mathcal{B}^{<k>}}$ & $:$ & $\mathcal{B}^{<k>}$ & $\rightarrow $ & $%
\mathcal{B}^{<k>},$ & $u$ & $\mapsto $ & $[\mathbf{n}\mapsto \mathrm{Rev}%
(u)(-\mathbf{n)]}.$%
\end{tabular}
\end{equation*}
Moreover $\widetilde{\mathcal{B}}^{<k>}\simeq R[\mathbf{x},\mathbf{x}%
^{-1}]^{\circ }$ becomes a Hopf $R$-algebra $(\widetilde{\mathcal{B}}%
^{<k>},\star _{g},\eta _{g},\Delta _{\widetilde{\mathcal{B}}%
^{<k>}},\varepsilon _{\widetilde{\mathcal{B}}^{<k>}},S_{\widetilde{\mathcal{B%
}}^{<k>}}),$ where $\star $ is the Hadamard product (\ref{Had-2}), $\eta $
is the unity (\ref{et-bi}) and 
\begin{equation*}
\begin{tabular}{llllllll}
$\Delta _{\widetilde{\mathcal{B}}^{<k>}}$ & $:$ & $\widetilde{\mathcal{B}}%
^{<k>}$ & $\rightarrow $ & $\widetilde{\mathcal{B}}^{<k>}\otimes _{R}%
\widetilde{\mathcal{B}}^{<k>},$ & $\widetilde{u}$ & $\mapsto $ & $%
\sum\limits_{\mathbf{t}\leq \mathbf{l}-\mathbf{1}}\mathrm{Rev}(\mathbf{x}^{%
\mathbf{t}}\rightharpoonup \widetilde{u}_{\mid _{\mathbb{N}%
_{0}^{k}}}))\otimes \mathrm{Rev}(e_{\mathbf{t}}^{\mathbf{F}}),$ \\ 
$\varepsilon _{\widetilde{\mathcal{B}}^{<k>}}$ & $:$ & $\widetilde{\mathcal{B%
}}^{<k>}$ & $\rightarrow $ & $R,$ & $\widetilde{u}$ & $\mapsto $ & $%
\widetilde{u}(\mathbf{0}),$ \\ 
$S_{\widetilde{\mathcal{B}}^{<k>}}$ & $:$ & $\widetilde{\mathcal{B}}^{<k>}$
& $\rightarrow $ & $\widetilde{\mathcal{B}}^{<k>},$ & $\widetilde{u}$ & $%
\mapsto $ & $[\mathbf{z}\mapsto \widetilde{u}(-\mathbf{z)]}.$%
\end{tabular}
\end{equation*}
Note that with these structures the isomorphism $\mathcal{B}^{<k>}\simeq 
\widetilde{\mathcal{B}}^{<k>}$ of Lemma \ref{bi=rev} turns to be an
isomorphism of Hopf $R$-algebras.
\end{punto}

The following theorem extends the corresponding result from the case of a
base field \cite[Page 124]{LT90} (see also \cite[14.15]{KKMN95}) to the case
of arbitrary \emph{artinian }ground rings:

\begin{theorem}
\label{hopf-rev}If $R$ is artinian, then there are isomorphisms of $R$%
-bialgebras 
\begin{equation}
{\normalsize \mathcal{L}}^{<k>}\simeq \mathcal{D}^{<k>}\oplus \mathbb{%
\widetilde{{\normalsize \mathcal{L}}}}^{<k>}=\mathcal{D}^{<k>}\oplus \mathbb{%
\widetilde{\mathcal{B}}}^{<k>}\simeq \mathcal{D}^{<k>}\oplus \mathcal{B}%
^{<k>}=\mathcal{D}^{<k>}\oplus \mathcal{R}^{<k>}.  \label{iso-bial}
\end{equation}
\end{theorem}

\begin{Beweis}
Consider the isomorphism $\mathcal{L}^{<k>}\simeq \mathcal{D}^{<k>}\oplus 
\widetilde{\mathcal{L}}^{<k>}$ (\ref{LM=}). With the help of Lemmata \ref
{q-q} and \ref{qx} one can show as in \cite[Seite 123]{LT90}, that $\gamma :%
\mathcal{L}^{<k>}\rightarrow \mathbb{\widetilde{\mathcal{L}}}^{<k>}$ (\ref
{gama}) and $\beta :\mathbb{\widetilde{\mathcal{L}}}^{<k>}\rightarrow 
\mathcal{L}^{<k>}$ (\ref{beta}) are in fact bialgebra morphisms. Obviously $%
\mathrm{Ke}(\gamma )=\mathcal{D}^{<k>}\subset \mathcal{L}^{<k>}$ is an $%
\mathcal{L}^{<k>}$-subbialgebra and we are done.$\blacksquare $
\end{Beweis}

As an analog to Corollary (\ref{LM*=}) we get

\begin{corollary}
Let $M$ be an $R[\mathbf{x},\mathbf{x}^{-1}]$-module. Then we have an
isomorphism of $R[\mathbf{x},\mathbf{x}^{-1}]^{\circ }$-comodules 
\begin{equation*}
\widetilde{\mathcal{L}}{\normalsize _{M^{\ast }}^{<k>}\simeq M[\mathbf{x}},%
\mathbf{x}^{-1}{\normalsize ]^{\circ }\simeq M^{\ast }\otimes _{R}R[\mathbf{x%
}},\mathbf{x}^{-1}{\normalsize ]^{\circ }\simeq M^{\ast }\otimes _{R}}\text{ 
}\widetilde{\mathcal{L}}_{R}^{<k>}{\normalsize .}
\end{equation*}
In particular $M[\mathbf{x},\mathbf{x}^{-1}]^{\circ }$ \emph{(}$\widetilde{%
\mathcal{L}}_{M^{\ast }}^{<k>}$\emph{)} is a cofree $R[\mathbf{x},\mathbf{x}%
^{-1}]^{\circ }$-comodule \emph{(}$\widetilde{\mathcal{L}}_{R}^{<k>}$%
-comodule\emph{)}.
\end{corollary}

As a consequence of \cite[Satz 2.4.7]{Abu2001} and \cite[Folgerung 2.5.10]
{Abu2001} we get

\begin{corollary}
\label{axk}Let $R$ be noetherian and consider the $R$-bialgebra $R[\mathbf{x}%
;g]^{\circ }$ \emph{(}resp. the Hopf $R$-algebra $R[\mathbf{x};p]^{\circ },$
the Hopf $R$-algebra $R[\mathbf{x},\mathbf{x}^{-1}]^{\circ }$\emph{).} If $A$
is an $\alpha $-algebra \emph{(}resp. an $\alpha $-bialgebra, a Hopf $\alpha 
$-algebra\emph{)}, then we have isomorphism of $R$-coalgebras \emph{(}resp. $%
R$-bialgebras, Hopf $R$-algebras\emph{)} 
\begin{equation}
A[\mathbf{x};g]^{\circ }\simeq A^{\circ }\otimes _{R}R[\mathbf{x};g]^{\circ
},\text{ }A[\mathbf{x};p]^{\circ }\simeq A^{\circ }\otimes _{R}R[\mathbf{x}%
;p]^{\circ }\text{ and }A[\mathbf{x},\mathbf{x}^{-1}]^{\circ }\simeq
A^{\circ }\otimes _{R}R[\mathbf{x},\mathbf{x}^{-1}]^{\circ }.
\label{A-Rx-iso}
\end{equation}
\end{corollary}

\begin{punto}
\textbf{Representative functions.} Let $G$ be a monoid (a group) and
consider the $R$-algebra $B=R^{G}$ with pointwise multiplication. Then $B$
is an $RG$-bimodule under the left and right actions 
\begin{equation*}
(yf)(x)=f(xy)\text{ and }(fy)(x)=f(yx)\text{ for all }x,y\in G.
\end{equation*}
We call $f\in R^{G}$ an $R$\emph{-valued representative function} on the
monoid $G,$ if $(RG)f(RG)$ is finitely generated as an $R$-module. If $R$ is
noetherian, then the subset $\mathcal{R}(G)\subset R^{G}$ of all
representative functions on $G$ is an $RG$-subbimodule. Moreover we deduce
from \cite[Theorem 2.13, Corollary 2.15]{AG-TW2000} that in case $%
(RG)^{\circ }\subset R^{G}$ is pure, we have an isomorphism of $R$%
-bialgebras (Hopf $R$-algebras) $\mathcal{R}(G)\simeq (RG)^{\circ }.$
\end{punto}

\begin{corollary}
Let $R$ be noetherian.

\begin{enumerate}
\item  Considering the monoid $(\mathbb{N}_{0}^{k},+)$ we have isomorphisms
of $R$-bialgebras 
\begin{equation*}
\mathcal{R}(\mathbb{N}_{0}^{k})\simeq R[\mathbf{x};p]^{\circ }\simeq 
\mathcal{L}_{R}^{<k>}.
\end{equation*}

\item  Considering the monoid $(\mathbb{Z}^{k},+)$ we have isomorphisms of
Hopf $R$-algebras 
\begin{equation*}
\mathcal{R}(\mathbb{Z}^{k})\simeq R[\mathbf{x},\mathbf{x}^{-1}]^{\circ
}\simeq \widetilde{\mathcal{B}}_{R}^{<k>}\simeq \mathcal{B}_{R}^{<k>}.
\end{equation*}
\end{enumerate}
\end{corollary}

\textbf{Acknowledgment.} Besides the new results, this note extends results
in Kapitel 4 of my doctoral thesis at the Heinrich-Heine Universit\"{a}t -
D\"{u}sseldorf (Germany). I am so grateful to my advisor Prof. Robert
Wisbauer for his wonderful supervision and the continuous encouragement and
support.

The author is also grateful to the referees for helpful suggestions and for
drawing his attention to the articles \cite{Kur00} and \cite{Kur94} by V.
Kurakin.

\end{document}